\documentclass[11pt,dvips]{article}
\usepackage{amsmath,amsfonts,latexsym,amssymb,chicago,eufrak,amscd, graphicx,pictexwd}
\newcommand{\mycite}[1]{{\small \sc \citeNP{#1}}}

\def\argmin{\mbox{argmin}}

\def\R{{\mathbb R}}

\def\F{{\cal F}}

\def\E{{\mathbb E}}

\def\labda1{\lambda_1}
\def\labda2{\lambda_2}

\def\argmin{\mbox{argmin}}

\def\comment#1{\relax}

\def\=in{\mathop{\rm =}}

\renewcommand{\P}{{\mathbb P}}
\renewcommand{\E}{{\mathbb E}}
\renewcommand{\argmin}{\mathop{\rm argmin}}

\renewcommand{\F}{{\mathbb F}}

\newcommand{\eps}{{\varepsilon }}

\newtheorem{thm}{Theorem}[section]
\newtheorem{coro}[thm]{Corollary}
\newtheorem{lem}[thm]{Lemma}
\newtheorem{prop}[thm]{Proposition}
\numberwithin{equation}{section}

\begin{document}
\title{Technical report: Adaptivity and optimality of the monotone least squares estimator for four different models.}
\author{Eric Cator\\ Delft University of Technology}
\date{\today}
\maketitle
\begin{abstract}
In this paper we will consider the estimation of a monotone regression (or density) function in a fixed point by the least squares (Grenander) estimator. We will show that this estimator is fully adaptive, in the sense that the attained rate is given by a functional relation using the underlying function $f_0$, and not by some smoothness parameter, and that this rate is optimal when considering the class of all monotone functions, in the sense that there exists a sequence of alternative monotone functions $f_1$, such that no other estimator can attain a better rate for both $f_0$ and $f_1$. We also show that under mild conditions the estimator attains the same rate in $L^q$ sense, and we give general conditions for which we can calculate a (non-standard) limiting distribution for the estimator.
\end{abstract}

\section{Introduction and results}

There exists an extensive literature on the problem of estimating a monotone increasing regression function or monotone decreasing density. We will consider the NPMLE or Grenander estimator for a monotone density, see \mycite{gren:56}, and the least squares estimator for a monotone regression function. Prakasa Rao obtained the rate and the limiting distribution for the Grenander estimator in a fixed point in \mycite{prakrao:69}, and in \mycite{brunk:70} a similar result was obtained for the least squares estimator. Results for global measures of convergence were obtained in \mycite{gretal:99} and \mycite{koulop:05} for the density case, and in \mycite{durot:05} for the regression case. A unified approach that incoorporates some other well known monotone estimators is given in \mycite{durot:07}. A common problem with these results is that they can only be proved under quite strong conditions, in which case there exist other non-isotonic estimators with faster rates.\\

\noindent Another approach, which addresses adaptivity, can be found in \mycite{kanglow:02}. Here the authors define an estimation procedure for $f_0(0)$, where $f_0$ is a monotone regression function in the white noise model, that is rate-adaptive in a minimax sense, for any $L_q$ $(q\geq 1)$ loss-function, with respect to a Lipschitz parameter $\alpha$. A serious drawback of this procedure compared to the estimator we consider, is that it does not, in general, give a monotone function as an estimate, when the procedure is applied to an interval of fixed points. Furthermore, the rate is described in terms of a Lipschitz parameter, not allowing for fast rates when the function $f_0$ has derivative $0$ in $0$, nor for slow rates when the function is not Lipschitz for any parameter value, nor for rates that cannot be described by the Lipschitz parameter alone (think of logarithmic corrections).\\

\noindent In this paper we will derive the rate at which the least squares estimator (or Grenander estimator) estimates $f_0(0)$, for a completely general monotone function $f_0$ that is continuous in $0$,  for four different models: the white noise model, measurements on a grid (not necessarily with normal errors), measurements on random design points and a sample from a decreasing density. This rate is defined in terms of the probabilistic error. This means that if we for example consider the white noise model
\[ Y(t) = \int_0^t f_0(s)\,ds + \frac1{\sqrt{n}} W(t),\]
for $t\in [-1,1]$, and we fix some $0<\alpha <1$, we get a rate $a_n$ that satisfies
\[ \P(|\hat{f}(0)-f_0(0)|\geq a_n) \leq \alpha.\]
This way of determining a rate is essential if we wish to get the full generality of our results, as was also observed in \mycite{cailow:06}. Our rate is defined in terms of a functional relation involving the function $f_0$. It turns out that the rate is similar for all four models, and to define it in the white noise case, we assume without loss of generality that $f_0(0)=0$ and define
\[ F_0(t) = \int_0^t f_0(s)\,ds.\]
Note that this function is convex, and due to the continuity of $f_0$ in $0$ (without which, we cannot estimate $f_0(0)$ consistently), we have that $F_0'(0)=0$. Then we fix $C>0$ and we define $a,r_a>0$ and $b,r_b>0$ depending on $n$ such that
\begin{equation}\label{eq:intro}
F_0(r_a)= ar_a,\ \ F_0(-r_b)= br_b\ \ \ \mbox{and}\ \ r_a^{1/2}a=r_b^{1/2}b=Cn^{-1/2}.
\end{equation}
Note that except for the simpler case where $F_0(-1)=0$ or $F_0(1)=0$, these equations always have a unique solution for $n$ large enough. Define the functions $\psi_l$ and $\psi_r$ by
\begin{equation}\label{eq:defpsi} \psi_r(s) = \limsup_{t\downarrow 0} \frac{F_0(st)}{F_0(t)} \ \ \mbox{and}\ \ \ \psi_l(s) = \limsup_{t\uparrow 0} \frac{F_0(st)}{F_0(t)}\ \ \ (s\in [0,1]).
\end{equation}
If $F_0(t)=0$ for some $t>0$, we define $\psi_r(s)=0$ for $s\in [0,1)$ and $\psi_r(1)=1$, and likewise for $\psi_l$. It is not hard to show that $\psi_r$ and $\psi_l$ are convex increasing functions, such that $0\leq \psi_r(s),\psi_l(s)\leq s$. We will show the following theorem:
\begin{thm}\label{thm:ratewn}
With the notations as above, we have that
\begin{eqnarray*}
\limsup_{n\to \infty}\ \P(\hat{f}(0)\geq a) & \leq  & \P\left(\inf_{s\leq 0} W(s) - Cs \leq \inf_{0\leq s\leq 1} W(s) - C(s-\psi_r(s))\right),\\
\limsup_{n\to \infty}\ \P(\hat{f}(0)\leq -b) & \leq  & \P\left(\inf_{s\leq 0} W(s) - Cs \leq \inf_{0\leq s\leq 1} W(s) - C(s-\psi_l(s))\right).
\end{eqnarray*}
\end{thm}
The actual rate is therefore given by $\max(a,b)$, since the probability on the right hand side always goes to zero as $C\to \infty$, and when $\psi_r$ or $\psi_l$ differs from the identity function, the respective probability goes to zero exponentially fast in $C$. This will also allow us to establish $L^q$ convergence of the estimator in the following sense. Define the increasing function
\[ G_0(t) = F_0(t)/t.\]
If $F_0(1)>0$, it is possible to define $G_0^{-1}$ as a strictly increasing continuous function on $[0,F_0(1)]$. Choose $a_\delta$ small enough and define
\[ H_0(a) = \left\{\begin{array}{ll}
a\sqrt{G_0^{-1}(a)} & \mbox{if } |a|\leq a_\delta\\
a-{\rm sgn}(a)a_\delta + H_0({\rm sgn}(a)a_\delta) & \mbox{if } |a|>a_\delta.\end{array}\right. \]
The connection with the rate equations \eqref{eq:intro} is given by $G_0^{-1}(a)=r_a$ and $H_0(a)=Cn^{-1/2}$, for $0<a<a_\delta$. We will show the following theorem:
\begin{thm}
Suppose $\psi_r(s)<s$ for some $s\in (0,1)$. Let $\chi:[0,\infty)\to [0,\infty)$ be such that for some constants $K>0$ and $m\geq 1$
\begin{eqnarray*} \chi(a) & \leq  & H_0(a) \ \ \ \ \ \hspace{0.5pt} \ \mbox{for } a\leq a_\delta \\
\chi(a) & \leq & KH_0(a)^m \ \ \mbox{for } a>a_\delta.
\end{eqnarray*}
Then there exists constants $L_1,L_2,\gamma, n_0>0$ such that for all $n\geq n_0$ and $C>0$
\[ \P(n^{-1/2}\chi(\hat{f}(0)_+)\geq C) \leq L_1e^{-L_2C^\gamma}.\]
\end{thm}
We can use this to show that if there exists $\alpha,M>0$ such that $f_0(x)\leq Mx^\alpha$ for positive $x$ in a neighborhood of $0$, then for any $q>0$
\[ \limsup_{n\to \infty}\ n^{\frac{\alpha}{2\alpha+1}}\,\E((\hat{f}(0)-f_0(0))_+^q) < +\infty.\]
Here we use the notation $x_+=\max(0,x)$. Note that controlling the behavior of $f_0$ to the right of $0$, only controls the ``overshoot'' of the estimator.\\

\noindent
We will also determine weak regularity conditions for $F_0$, such that we can determine the limiting distribution of $\hat{f}(0)$. Suppose
\[ \lim_{n\to \infty} \frac{r_a}{r_b} = \gamma \in [0,\infty).\]
This says that the rates to the left and to the right of $0$ are well behaved with respect to each other, which is a natural condition for a limiting distribution to exist. Furthermore, suppose $F_0$ is regularly varying near $0$: there exists $\alpha > 1$, such that for all $s>0$
\[ \lim_{t\downarrow 0} \frac{F_0(st)}{F_0(t)} = s^\alpha.\]
This says that $F_0$ scales properly near $0$, which is another natural condition: we don't want different behavior of $F_0$ for different scales. We will prove the following theorem:
\begin{thm}
Let $W_s$ $(s\in \R)$ denote twosided standard Brownian motion, and define the process
\[ X(s) = \left\{ \begin{array}{ll}
W_s + s^\alpha & \mbox{for } s\geq 0,\\
W_s + \gamma^{\alpha-1/2}|s|^\alpha & \mbox{for } s\leq 0, \end{array}\right.\]
and the process $\hat{X}(s)$ as the greatest convex minorant of $X$. With the conditions given above, we have that
\[ \frac{\hat{f}(0)_+}{H_0^{-1}(n^{-1/2})} \stackrel{d}{\longrightarrow} \frac{d\hat{X}}{ds}(0)_+ \ \ \mbox{and}\ \  \frac{\hat{f}(0)_-}{-H_0^{-1}(-n^{-1/2})} \stackrel{d}{\longrightarrow} \frac{d\hat{X}}{ds}(0)_-.\]
\end{thm}

\noindent
Finally we will show that the rate for $\hat{f}(0)$ is local asymptotic minimax:
\begin{thm}\label{thm:optimality}
Choose two significance levels $\alpha \in (0,1)$ and $\beta \in (0,1/2)$. There exist $\eta>0$, such that for all $n$ large enough, we can find a monotone function $f_1$ (close to $f_0$), and we can find a rate $\gamma_n$ with
\[ \limsup_{n \to \infty}\ \max_{i=0,1}\ \P_{f_i}\left(|\hat{f}(0)-f_i(0)|\geq \gamma_n\right) \leq \alpha\]
and
\[  \liminf_{n \to \infty}\ \inf_{\hat{\theta}}\ \max_{i=0,1}\ \P_{f_i}\left(|\hat{\theta}(Y)-f_i(0)|\geq \eta\cdot \gamma_n\right) > \beta,\]
where $\hat{\theta}(Y)$ is any estimator of $f(0)$ based on the data $Y$.
\end{thm}
This says that $\hat{f}(0)$ attains a certain rate $\gamma_n$ for both $f_0$ and the sequence of alternatives $f_1$ (of course we take $\gamma_n$ a constant times $\max(a,b)$), and no other estimator can do significantly better for both $f_0$ and $f_1$ simultaneously. This way of describing optimality was inspired by a talk in Oberwolfach, given by Tony Cai and Mark Low, although their concept looked at the $L_q$-risk and it did not require that $\hat{f}(0)$ estimate $f_1$ with the same rate.\\

\noindent We would like to give some feel for Equations \eqref{eq:intro}. Suppose $f_0$ is Lipschitz continuous in $0$ with parameter $\alpha >0$, so for $x$ in a neighbourhood of $0$, we have (remember that $f_0(0)=0$)
\[ |f_0(x)|\lesssim |x|^\alpha.\]
Here, $g(x) \lesssim h(x)$ denotes there exists a constant $M>0$ such that $g(x)\leq Mh(x)$ for all relevant $x$. Then $F_0(x)\lesssim |x|^{\alpha +1}$, so \eqref{eq:intro} gives us
\[ ar_a\lesssim r_a^{\alpha + 1}.\]
This means that $r_a^{-1}\lesssim a^{-1/\alpha}$. Together with the second equality for $a$ in \eqref{eq:intro}
\[ a \lesssim r_a^{-1/2}n^{-1/2},\]
this leads to
\[ a\lesssim n^{-\frac{\alpha}{2\alpha +1}}.\]
For $b$ we can derive the same bound. This corresponds to the rate found in \mycite{kanglow:02}. Another interesting case is when $\lim_{a\to 0} r_a=r_0>0$. This means that $f_0$ is flat to the right of $0$, on the interval $[0,r_0)$. Then
\[ a \lesssim r_0^{-1/2}n^{-1/2},\]
so this corresponds to a parametric rate.\\

\noindent The rest of the paper is organized as follows: Section \ref{sec:psi} considers the functions $\psi_r$ and $\psi_l$. Section \ref{sec:wn} deals with the white noise model, for which we will prove all the above results. Sections \ref{sec:grid}, \ref{sec:rp} and \ref{sec:gren} each deal with one of the other three models we will consider, but we will only formulate and prove the corresponding theorems of Theorem \ref{thm:ratewn}, giving the rate, and Theorem \ref{thm:optimality}, showing the optimality of the rate. The other theorems, concerning the $L^q$ convergence and the limiting distribution, require some weak technical conditions, but the ideas are the same as for the white noise model, and are not worked out in this paper.

\section{The functions $\psi_r$ and $\psi_l$}\label{sec:psi}

In this section we will take a closer look to the functions $\psi_r$ and $\psi_l$ defined in \eqref{eq:defpsi} in the Introduction. We will concentrate on $\psi_r$, since completely analogous statements will hold for $\psi_l$. Since the function $s\mapsto F_0(st)$ is convex and increasing for all $t>0$, we get that $\psi_r(s)$ is also an increasing and convex function on $[0,1]$ (this is true for the $\limsup$ of convex functions, not necessarily for the $\liminf$). Furthermore, we clearly have that $\psi_r(0)=0$ and $\psi_r(1)=1$. Finally, since $F_0$ is convex, we know that for $s\in [0,1]$,
\[ F_0(st)\leq sF_0(t) + (1-s)F_0(0) = sF_0(t).\]
This shows that for any $F_0$ we have that $\psi_r(s)\leq s$.
\begin{lem}\label{lem:unifconvpsi}
For each $\tau\in [0,1)$, there exists a positive continuous increasing function $\eta$ with
$$F_0(t)=0\Rightarrow \eta(t)=0\ \ \ (\forall t\in [0,1]),$$
such that for all $t\in (0,1]$ and for all $s\in [0,\tau]$
\[\ F_0(st)\leq (\psi_r(s) + \eta(t))F_0(t).\]
\end{lem}
{\bf Proof:} Suppose $F_0(t)>0$ for all $t>0$. Define the auxiliary functions
\[ G_t(s) = \sup_{u\leq t} \frac{F_0(su)}{F_0(u)}.\]
These functions are all convex and they decrease pointwise to $\psi_r$ on $[0,1]$. Since $G_t(0)=0$ for all $t>0$, we conclude that $G_t$ converges uniformly to $\psi_r$ on $[0,\tau]$. Define
\[ \eta(0)=0\ \ \ \mbox{and}\ \ \ \eta(t)=\sup_{s\in [0,\tau]} |\psi_r(s) - G_t(s)|\ \ (t\in(0,1])\]
and note that
\[ \frac{F_0(st)}{F_0(t)} \leq G_t(s) \leq \psi_r(s) + \eta(t),\]
to conclude the statement of the lemma (note that $\psi_r$ is continuous on $[0,\tau]$, so $\eta$ is indeed a continuous increasing function). Now suppose that $F_0(t)=0$ for some $t>0$. We defined $\psi_r(s)=0$ for $s\in [0,1)$ in this case. Define
\[ r_0 = \sup\{ t\in [0,1] : F_0(t)=0\}.\]
If $r_0=1$, then the statement of the lemma holds with $\eta=0$. Suppose $r_0<1$. Since $s\leq \tau$, we have for each $t\leq r_0/\tau$, $F_0(st)=0$. Define $\eta(t)=0$ for $t\in [0,r_0]$, $\eta(t)=1$ for $t\in [r_0/\tau , 1]$, and continuous in between, and the statement of the lemma holds trivially.\hfill $\Box$\\

Let $\{W_s\,:\,s\in \R\}$ be a two-sided Brownian motion. We will encounter the probabilities in the next lemma throughout the rest of the paper.
\begin{lem}\label{lem:WpsiC}
For any $F_0$ we have that
\[ \P\left(\inf_{s\leq 0} W_s - Cs \leq \inf_{0\leq s\leq 1} W_s - C(s-\psi_r(s))\right) \leq \frac1{\sqrt{2\pi}\,C}.\]
If there exists $s\in (0,1)$ such that $\psi_r(s)<s$, then there exist $\tau \in (0,1)$ and $\rho \in (0,1]$ such that
\[ \P\left(\inf_{s\leq 0} W_s - Cs \leq \inf_{0\leq s\leq 1} W_s - C(s-\psi_r(s))\right) \leq \sqrt{\frac{2}{\pi\tau}}\,\frac{1}{C\rho(2-\rho)}\,e^{-C^2\tau \rho^2/2}.\]
\end{lem}
{\bf Proof:} First note that the left-hand side and the right-hand side of two-sided Brownian motion are independent. It is therefore enough to consider the two sides within the probability seperately. It is well known that
\[ \P(\inf_{s\leq 0} W_s - Cs\leq -v) = \P(\sup_{s\geq 0} W_s - Cs \geq v) = e^{-2Cv}.\]
This follows from the hitting time of a linear boundary. Since for all $F_0$ we have that $\psi_r(s)\leq s$, we also need that
\[ \P(\inf_{0\leq s\leq 1} W_s \leq -w) = \P(\sup_{0\leq s\leq 1} W_s \geq w) = 2(1-\Phi(w)),\]
where $\Phi$ is the distribution function of the standard normal distribution. We get
\begin{eqnarray*}
\P\left(\inf_{s\leq 0} W_s - Cs \leq \inf_{0\leq s\leq 1} W_s - C(s-\psi_r(s))\right) & \leq & \P\left(\inf_{s\leq 0} W_s - Cs \leq \inf_{0\leq s\leq 1} W_s\right) \\
& = & \frac{2}{\sqrt{2\pi}}\,\int_0^\infty e^{-2Cw}e^{-\frac12 w^2}\,dw \\
& = & 2e^{2C^2}(1-\Phi(2C))\\
& \leq & \frac{1}{\sqrt{2\pi}\,C}.
\end{eqnarray*}
Now suppose that for some $s\in (0,1)$, $\psi_r(s)<s$. Since $\psi_r$ is convex, $\psi_r(0)=0$ and $\psi_r(1)=1$, this implies that for any $\tau\in (0,1)$ and any $s\in (0,\tau]$, $\psi_r(s)\leq s\psi_r(\tau)/\tau<s$. Choose $\tau \in (0,1)$ and define $\rho = 1 - \psi_r(\tau)/\tau>0$. Then
\[ \forall s\in [0,\tau]:\ s-\psi_r(s)\geq \rho s.\]
Now use that
\begin{eqnarray*}
\P\left(\inf_{s\leq 0} W_s - Cs \leq \inf_{0\leq s\leq 1} W_s - C(s-\psi_r(s))\right) & \leq & \P\left(\inf_{s\leq 0} W_s - Cs \leq \inf_{0\leq s\leq \tau} W_s - C\rho s\right) \\
& \leq & \P\left(\inf_{s\leq 0} W_s - Cs \leq W_\tau - C\rho \tau\right).
\end{eqnarray*}
This last probability we can calculate exactly:
\begin{eqnarray}\label{eq:ineqC}
\P\left(\inf_{s\leq 0} W_s - Cs \leq W_\tau - C\rho \tau\right) & = & 1-\Phi(C \sqrt{\tau} \rho) + \frac1{\sqrt{2\pi \tau}}\,\int_{0}^\infty e^{-2Cv}e^{-\frac{(v- C\tau \rho )^2}{2\tau}}\,dv\\ \nonumber
& = & 1-\Phi(C \sqrt{\tau} \rho ) + (1-\Phi(C \sqrt{\tau}(2-\rho) ))\,e^{C^2\tau(2-\rho)^2/2}\,e^{-C^2\tau\rho^2/2}\\ \nonumber
& \leq & \sqrt{\frac{2}{\pi\tau}}\,\frac{1}{C\rho(2-\rho)}\,e^{-C^2\tau \rho^2/2}.\hspace{6cm} \Box
\end{eqnarray}

\noindent
We will now consider the case where $F_0(st)/F_0(t)$ actually has a limit. This is comparable to saying that $F_0$ is a regularly varying function in $0$, but we have the extra information that $F_0$ is convex.
\begin{lem}\label{lem:limpsi}
Suppose for each $s\in (0,1]$ we have $F_0(s)>0$ and
\[ \lim_{t\downarrow 0} \frac{F_0(st)}{F_0(t)} = \psi_r(s).\]
Then either $\psi_r(s)=0$ on $[0,1)$, or $\psi_r(s)=s^\alpha$, for some $\alpha \geq 1$. In the latter case, we have that for each $\tau > 0$ (also for $\tau\geq 1$)
\[ \sup_{s\in [0,\tau]} \left( \frac{F_0(st)}{F_0(t)} - s^\alpha\right) \stackrel{t\downarrow 0}{\longrightarrow} 0.\]
\end{lem}
{\bf Proof:} Suppose $0\leq u\leq s\leq 1$. Then
\[ \psi_r(us)=\lim_{t\downarrow 0} \frac{F_0(ust)}{F_0(t)} = \lim_{t\downarrow 0} \frac{F_0(ust)}{F_0(st)}\frac{F_0(st)}{F_0(t)} = \psi_r(u)\psi_r(s).\]
Since $\psi_r$ is continuous and convex on $[0,1)$ and $\psi_r(0)=0$, we conclude that either $\psi_r(s)=0$ on $[0,1)$ or $\psi_r(s) = s^\alpha$ with $\alpha\geq 1$. In this last case, choose $s>1$. Then
\[ \lim_{t\downarrow 0} \frac{F_0(st)}{F_0(t)}= \left(\lim_{t\downarrow 0} \frac{F_0(t)}{F_0(st)}\right)^{-1} = \left(\lim_{t\downarrow 0} \frac{F_0(s^{-1}t)}{F_0(t)}\right)^{-1} = s^\alpha.\]
The family of convex functions $\{s\mapsto F_0(st)/F_0(t)\}$ converges pointwise to the convex function $s\mapsto s^\alpha$, and all functions are 0 in 0, so the convergence is actually uniform on compact subsets of $[0,\infty)$.\hfill $\Box$\\

\section{The monotone LS-estimator in white noise}\label{sec:wn}
\label{sec:heurrate}

We will work in the white noise model, so our data $Y(t)$ satisfies
\[ dY(t) = f_0(t)dt + \eps dW(t),\]
where $f_0$ is a monotone $L^2$-function on $[-1,1]$ and $W(t)$ is standard two-sided Brownian motion. As usual, the parameter $\eps$ should be compared to $n^{-1/2}$. We wish to study the least squares estimator, but in fact we will define for a realization of $W(t)$,
\[ Y(t) = \int_0^t f_0(t)dt + \eps W(t),\]
and the convex function
\[ \hat{F}(t) = \sup \{ \phi(t) : \phi \mbox{ affine and } \forall s\in [-1,1]: \phi(s)\leq Y(s)\}.\]
So $\hat{F}$ is the greatest convex minorant of $Y$. Now we define the estimator $\hat{f}$ as the left-derivative of the convex function $\hat{F}$, so for $t\in (-1,1)$
\[ \hat{f}(t) = \lim _{h\downarrow 0} \frac{\hat{F}(t)-\hat{F}(t-h)}{h}.\]
This is a monotone function and can be seen as a limit of least squares estimators over the class of monotone functions absolutely bounded by $M$, as $M\to \infty$.

We will assume without loss of generality that $f_0(0)=0$. Furthermore, to ensure that our estimator $\hat{f}(0)$ is consistent as $\eps \to 0$, we assume that $f_0$ is continuous in $0$. We are interested in the probability of the event $\{\hat{f}(0)\geq a\}$, for $a>0$. Define
\[ F_0(t) = \int_0^t f_0(s)ds.\]
Fix $C>0$ not depending on $\eps$, and choose $a,b>0$ and $r_a,r_b>0$ such that
\begin{equation}\label{eq:rate}
F_0(r_a)= ar_a,\ \ F_0(-r_b)= br_b\ \ \mbox{and}\ \ r_a^{1/2}a=r_b^{1/2}b=C\eps.
\end{equation}
Since $F_0$ is convex and continuous, and $f_0$ is continuous in $0$, this can always be done if $F_0(1)>0$ and $F_0(-1)>0$, simply by choosing $\eps$ small enough. We will consider the special (and simpler) case $f_0(t)=0$ for all $t>0$ (or for all $t<0$) separately.
\begin{thm}\label{thm:wnuppbnd}
With the notations as above, we have that
\[ \limsup_{\eps\downarrow 0} \P(\hat{f}(0)\geq a) \leq \P\left(\inf_{s\leq 0} W_{s}- Cs \leq  \inf_{0< s\leq 1} W_{s} - C(s-\psi_r(s))\right)\]
and
\[ \limsup_{\eps\downarrow 0} \P(\hat{f}(0)\leq -b) \leq \P\left(\inf_{s\leq 0} W_{s}- Cs \leq  \inf_{0< s\leq 1} W_{s} - C(s-\psi_l(s))\right).\]
Since both probabilities tend to zero when $C\to \infty$, it follows that Equations (\ref{eq:rate}) determine an upper bound for the rate of convergence of $\hat{f}(0)$.
\end{thm}
{\bf Proof:} We will only show the result for $a$; the proof for $b$ is completely similar. Note that we have the following ``switch relation'' for the greatest convex minorant:
\begin{equation}\label{eq:max<0}
\{\hat{f}(0)\geq a\} = \{ \inf_{-1\leq t\leq 0}(\eps W_t+F_0(t)-at) \leq \inf_{0< t\leq 1}(\eps W_t+F_0(t)-at)\}.
\end{equation}
We can rewrite (\ref{eq:max<0}) as follows:
\begin{eqnarray}\label{eq:rewrate}
\{\hat{f}(0)\geq a\} & = & \{ \inf_{-r_a^{-1}\leq s\leq 0}(r_a^{-1/2}W_{r_a s} + \eps^{-1}r_a^{-1/2}F_0(r_a s) - \eps^{-1}r_a^{1/2}as) \leq \nonumber \\
&  & \hspace{2cm} \inf_{0< s\leq r_a^{-1}}(r_a^{-1/2}W_{r_a s} + \eps^{-1}r_a^{-1/2}F_0(r_a s) - \eps^{-1}r_a^{1/2}as)\}.
\end{eqnarray}
Define
\[ \tilde{W}_{s} = r_a^{-1/2}W_{r_a s}.\]
Clearly, $\tilde{W}_s$ is also a two-sided Brownian motion. Now we can use Lemma \ref{lem:unifconvpsi}: for any $\tau\in(0,1)$ there exists a positive continuous function $\eta$ with $F_0(t)=0 \Rightarrow \eta(t)=0$, such that
\[ \inf_{0< s\leq r_a^{-1}}(r_a^{-1/2}W_{r_a s} + \eps^{-1}r_a^{-1/2}F_0(r_a s) - \eps^{-1}r_a^{1/2}as) \leq \inf_{0< s\leq \tau} \tilde{W}_{s} - C(s-\psi_r(s)) + C\eta(r_a).\]
Now remark that for $s<0$, $F_0(s)\geq 0$, so that
\[\inf_{-r_a^{-1}\leq s\leq 0}(r_a^{-1/2}W_{r_a s} + \eps^{-1}r_a^{-1/2}F_0(r_a s) - Cs) \geq \inf_{s\leq 0}(\tilde{W}_{s}- Cs).\]
In view of \eqref{eq:rewrate}, we have shown that
\begin{eqnarray}\label{eq:probW}
\P(\hat{f}(0)\geq a) & \leq & \P\left(\inf_{s\leq 0}(\tilde{W}_{s}- Cs) \leq  \inf_{0< s\leq \tau} \tilde{W}_{s} - C(s-\psi_r(s)) + C\eta(r_a)\right).
\end{eqnarray}
Define $r_0 = \lim_{a\downarrow 0} r_a$. Since we always have that $F_0(r_0)=0$, we conclude that $\lim_{a\downarrow 0} \eta(r_a) = 0$, so
\[ \limsup_{\eps \downarrow 0} \P(\hat{f}(0)\geq a) \leq \P\left(\inf_{s\leq 0}(\tilde{W}_{s}- Cs) \leq  \inf_{0< s\leq \tau} \tilde{W}_{s} - C(s-\psi_r(s))\right).\]
Since this is true for all $\tau \in (0,1)$, and since $\psi_r$ is increasing on $[0,1]$, we conclude that
\[ \limsup_{\eps \downarrow 0} \P(\hat{f}(0)\geq a) \leq \P\left(\inf_{s\leq 0}(\tilde{W}_{s}- Cs) \leq  \inf_{0< s\leq 1} \tilde{W}_{s} - C(s-\psi_r(s))\right).\]

When $f_0(t)=0$ for all $t>0$, we choose $a=C\eps$, and (\ref{eq:max<0}) implies that
\begin{eqnarray*}
 \P\left(\hat{f}(0)\geq a\right) & \leq & \P\left(\inf_{-1\leq t\leq 0}(W_t-Ct) \leq \inf_{0< t\leq 1}(W_t-Ct)\right).
 \end{eqnarray*}
This shows that in this case, the upper confident limit for $\hat{f}(0)$ is of order $\eps$ (parametric rate). This also happens when $r_0>0$, which is the case when $f_0$ is flat to the right of $0$.\hfill $\Box$\\

\subsection{$L^q$ convergence of the LS-estimator}

The basis for deriving the $L^q$ ($q>0$) convergence of the Least Squares estimator will be Equation \eqref{eq:probW}, together with Lemma \ref{lem:WpsiC} and a uniform integrability argument. We note that \eqref{eq:probW} holds for all choices of $C>0$, as long as Equations \eqref{eq:rate} for $a$ and $r_a$ have a solution. To ensure this, we choose $\delta\in (0,1)$ small and define $a_\delta>0$ and $C_\delta>0$ such that
\[ F_0(\delta)=a_\delta \delta\ \ \ \mbox{and}\ \ \ a_\delta \delta^{1/2} = \eps C_\delta.\]
This is possible as soon as $F_0(1)>0$; the case $F_0(1)=0$ is in fact easier. So for any $C\leq C_\delta$, we have $r_a\leq \delta$, and since $\eta$ is increasing, we get
\begin{eqnarray*}
\P(\hat{f}(0)\geq a) & \leq & \P\left(\inf_{s\leq 0}(W_{s}- Cs) \leq  \inf_{0< s\leq \tau} W_{s} - C(s-\psi_r(s)) + C\eta(\delta)\right).
\end{eqnarray*}
Note that from the derivation of Equation \eqref{eq:probW}, it follows that we can choose $\tau\in (0,1)$ fixed, independent of $C$ and $\eps$. Now we wish to make the dependence of $a$ on $C$ more specific. To this end, we introduce two auxiliary functions $G_0$ and $H_0$:
\[ G_0(t) = F_0(t)/t\ \ \ \ \mbox{and}\ \ \ \ H_0(a) = a\sqrt{G_0^{-1}(a)}\ \ \ \mbox{for } a\leq a_\delta.\]
Since $F_0$ is convex and has derivative $0$ in $0$, we get that $G_0$ is strictly increasing on the set $\{t\in [0,1] : F_0(t)>0\}$. This means that $G_0^{-1}$ can be defined on $[0,a_\delta]$ in a continuous, strictly increasing manner. Therefore, $H_0$ will be a continuous, strictly increasing function on $[0,a_\delta]$ with $H_0(0)=0$. Clearly,
\[r_a=G_0^{-1}(a)\ \ \ \mbox{and}\ \ \  H_0(a)=C\eps.\]
We extend the definition of $H_0$ to $[0,\infty)$ by
\[ H_0(a) = H_0(a_\delta) + a - a_\delta\ \ \ \mbox{for } a\geq a_\delta.\]
In this way, $H_0$ remains a continuous and strictly increasing function. We define for all $C>0$ and $\eps >0$, $H_0(a) = C\eps$. We can show the following proposition, using the notation $x_+=\max(0,x)$.
\begin{prop}\label{prop:H0f}
Suppose $\psi_r(s)<s$ for some $s\in (0,1)$ (and hence for all $s\in (0,1)$). With the notations as above, we have for all $\eps$ small enough, that for all $C>0$
\[ \P(\eps^{-1}H_0(\hat{f}(0)_+)\geq C) \leq \P\left(\inf_{s\leq 0}(W_{s}- Cs) \leq W_{\tau} - \frac14 \delta C\right).\]
\end{prop}
{\bf Proof:} Since $H_0$ is strictly increasing, we have
\[ \P(\hat{f}(0)\geq a) = \P(H_0(\hat{f}(0)_+) \geq H_0(a))=\P(\eps^{-1}H_0(\hat{f}(0)_+)\geq C).\]
So for any $C\leq C_\delta = a_\delta \delta^{1/2}\eps^{-1}$, we get
\[ \P(\eps^{-1}H_0(\hat{f}(0)_+)\geq C) \leq \P\left(\inf_{s\leq 0}(W_{s}- Cs) \leq \inf_{0< s\leq \tau} W_{s} - C(s-\psi_r(s)) + C\eta(\delta)\right).\]
What can we say when $C>C_\delta$? With $a$ defined by $H_0(a)=C\eps$, we get that $a>a_\delta$ and
\[ \P(\eps^{-1}H_0(\hat{f}(0))_+\geq C) = \P(\hat{f}(0)\geq a),\]
so we can use Equation \eqref{eq:max<0} to conclude that
\begin{eqnarray*}
\P(\eps^{-1}H_0(\hat{f}(0)_+)\geq C) & = & \P \left( \inf_{-1\leq t\leq 0}(\eps W_t+F_0(t)-at) \leq \inf_{0< t\leq 1}(\eps W_t+F_0(t)-at)\right)\\
& \leq & \P \left( \inf_{-1\leq t\leq 0}(\eps W_t-at) \leq \inf_{0< t\leq \delta}(\eps W_t+a_\delta t-at)\right)\\
& \leq & \P \left( \inf_{-1\leq t\leq 0}(W_t-(a_\delta \eps^{-1} + C-C_\delta)t) \leq \inf_{0< t\leq \delta}(W_t - (C-C_\delta)t)\right).
\end{eqnarray*}
Now note that $a_\delta \eps^{-1} = \delta^{-1/2}C_\delta > C_\delta$, which leads us to
\[ \P(\eps^{-1}H_0(\hat{f}(0)_+)\geq C) \leq  \P \left( \inf_{-1\leq t\leq 0}(W_t-Ct) \leq \inf_{0< t\leq \delta}(W_t - (C-C_\delta)t)\right).\]
Now we use that $\psi_r(s)<s$, also in view of Lemma \ref{lem:WpsiC}. Choose $\delta$ so small, that $\delta < \tau$ and $\tau - \psi_r(\tau)-\eta(\delta) > \delta/2$. Then for $C\leq C_\delta$ we get
\begin{eqnarray*}
 \P(\eps^{-1}H_0(\hat{f}(0)_+)\geq C) & \leq & \P\left(\inf_{s\leq 0}(W_{s}- Cs) \leq \inf_{0< s\leq \tau} W_{s} - C(s-\psi_r(s)) + C\eta(\delta)\right)\\
 & \leq & \P\left(\inf_{s\leq 0}(W_{s}- Cs) \leq W_{\tau} - \frac12 \delta C\right).
\end{eqnarray*}
For $C>C_\delta$ we have
\[ \P(\eps^{-1}H_0(\hat{f}(0)_+)\geq C) \leq \P \left( \inf_{t\leq 0}(W_t-Ct) \leq W_\delta - (C-C_\delta)\delta\right).\]
Since for $C$ big enough, we have
\[ \P\left(\inf_{s\leq 0}(W_{s}- Cs) \leq W_{\delta} - \frac12 \delta C\right) \leq \P\left(\inf_{s\leq 0}(W_{s}- Cs) \leq W_{\tau} - \frac12 \delta C\right),\]
we conclude that for $C\geq 2C_\delta$,
\[ \P(\eps^{-1}H_0(\hat{f}(0)_+)\geq C) \leq \P\left(\inf_{s\leq 0}(W_{s}- Cs) \leq W_{\tau} - \frac12 \delta C\right).\]
Note that $\P(\eps^{-1}H_0(\hat{f}(0)_+)\geq C)$ is a decreasing function of $C$, so for $\eps$ small enough, which means $C_\delta$ big enough, we can conclude for all $C>0$ that
\[ \P(\eps^{-1}H_0(\hat{f}(0)_+)\geq C) \leq \P\left(\inf_{s\leq 0}(W_{s}- Cs) \leq W_{\tau} - \frac14 \delta C\right).\]
\hfill $\Box$\\

The definition of $H_0(a)$ depends on the choice of $\delta>0$, when $a>a_\delta$, but since $\hat{f}(0)\to 0$, this should not be relevant. To prove this, we show the following corollary.
\begin{coro}\label{corr:chiH0}
Suppose $\psi_r(s)<s$ for some $s\in (0,1)$. Let $\chi:[0,\infty)\to [0,\infty)$ be such that for some constants $K>0$ and $n\geq 1$
\begin{eqnarray*} \chi(a) & \leq  & H_0(a) \ \ \ \ \ \hspace{0.5pt} \ \mbox{for } a\leq a_\delta \\
\chi(a) & \leq & KH_0(a)^n \ \ \mbox{for } a>a_\delta.
\end{eqnarray*}
Then there exists constants $L_1,L_2,\gamma, \eps_0>0$ such that for all $0<\eps<\eps_0$ and $C>0$
\[ \P(\eps^{-1}\chi(\hat{f}(0)_+)\geq C) \leq L_1e^{-L_2C^\gamma}.\]
\end{coro}
{\bf Proof:} Note that for $\eps < 1$,
\begin{eqnarray*}
\P(\eps^{-1}\chi(\hat{f}(0)_+)\geq C) & = & \P(\eps^{-1}\chi(\hat{f}(0)_+)\geq C \wedge \hat{f}(0)\leq a_\delta) + \P(\eps^{-1}\chi(\hat{f}(0)_+)\geq C\wedge \hat{f}(0)> a_\delta)\\
& \leq & \P(\eps^{-1}H_0(\hat{f}(0)_+)\geq C) + \P(K\eps^{-1}H_0(\hat{f}(0)_+)^n\geq C)\\
& \leq & \P(\eps^{-1}H_0(\hat{f}(0)_+)\geq C) + \P(\eps^{-1/n}H_0(\hat{f}(0)_+)\geq K^{-1/n}C^{1/n})\\
& \leq & \P(\eps^{-1}H_0(\hat{f}(0)_+)\geq C) + \P(\eps^{-1}H_0(\hat{f}(0)_+)\geq K^{-1/n}C^{1/n})\\
&\leq & L_1e^{-L_2C^\gamma},
\end{eqnarray*}
for some choice of $L_1, L_2, \gamma ,\eps_0>0$ and all $\eps<\eps_0$. In the last step we used Proposition \ref{prop:H0f} and Equation \eqref{eq:ineqC}.\hfill $\Box$\\

To see how we can use Corollary \ref{corr:chiH0}, let us assume that $\psi_r(s)<s$ and that for positive $x$ in  some neighbourhood of $0$, we have
\[ f_0(x) \lesssim x^\alpha,\]
for some $\alpha>0$. Then $F_0(x)\lesssim x^{\alpha +1}$, so $G_0(x)\lesssim x^\alpha$. Therefore, $G^{-1}(a)\gtrsim a^{1/\alpha}$, and
\[ H_0(a)\geq  Ra^{\frac{2\alpha + 1}{2\alpha}},\]
for some $R>0$, at least for $0\leq a\leq a_\delta$ if we choose $\delta>0$ small enough. Now define for all $a>0$
\[ \chi(a) = Ra^{\frac{2\alpha + 1}{2\alpha}}.\]
Corollary \ref{corr:chiH0} then shows that for all $C>0$
\[
\P(\eps^{-1}\chi(\hat{f}(0)_+) \geq C)  = \P\left( \eps^{-\frac{2\alpha}{2\alpha + 1}}\hat{f}(0) \geq (C/R)^{\frac{2\alpha}{2\alpha + 1}}\right) \leq L_1e^{-L_2C^\gamma},\]
and this proves that for any $q>0$,
\[ \lim_{\eps \to 0} \E\left(\left(\eps^{-\frac{2\alpha}{2\alpha + 1}}\hat{f}(0)_+\right)^q\right) < +\infty.\]
Of course we can get a similar result for $\hat{f}(0)_-$ and for $|\hat{f}(0)|$. The condition $\psi_r(s)<s$ does not affect the rate of the estimator. It is only necessary to get control over the tail of the rescaled LS estimator.

\subsection{Limiting distribution of the Least Squares estimator}

Our methods also allow us to derive non-standard limiting distributions for the Least Squares estimator. These limiting distributions only exist when $f_0$ is somehow ``regular'' near $0$. The precise conditions are described in the following theorem and will use Lemma \ref{lem:limpsi}. We start with the rate equations: for $\eps>0$ and $C>0$ we define $a, r_a, b$ and $r_b$ by
\[F_0(r_a)= ar_a,\ \ F_0(-r_b)= br_b\ \ \mbox{and}\ \ r_a^{1/2}a=r_b^{1/2}b=C\eps.\]

\begin{thm}\label{thm:limdist}
Suppose that
\[ \lim_{\eps\to 0} \frac{r_a}{r_b} = \gamma,\]
with $\gamma \in [0,\infty)$. Furthermore, suppose that for $s\geq 0$,
\[ \lim_{t\downarrow 0} \frac{F_0(st)}{F_0(t)} = s^{\alpha}\]
for $\alpha> 1$ (see also Lemma \ref{lem:limpsi}). Then, if $W_s$ $(s\in \R)$ denotes twosided standard Brownian motion,
\[ \lim_{\eps \to 0} \P(\hat{f}(0)>a)  =  \P \left( \inf_{s\leq 0}( W_s + C\gamma^{\alpha-1/2} |s|^{\alpha} - Cs ) \leq \inf_{s\geq 0} ( W_s + C |s|^{\alpha} - Cs ) \right)\]
If $\lim_{\eps\to 0} r_a/r_b = +\infty$, then
\[ \lim_{\eps \to 0} \P(\hat{f}(0)> 0)  = 0.\]
\end{thm}
{\bf Proof:} We start with assuming that $\gamma > 0$. Since $r_a/r_b \to \gamma$ and $ar_a^{1/2}=br_b^{1/2}$, we see that $a/b\to \gamma^{-1/2}$ and $F_0(r_a)/F_0(-r_b) \to \gamma^{1/2}$ (since $F_0(r_a)=ar_a$ and $F_0(-r_b)=br_b$). For each $\eta>0$, we have that $(\gamma - \eta)r_b \leq r_a\leq (\gamma + \eta)r_b$, for $\eps$ small enough. Therefore
\[ \limsup_{\eps \to 0} \frac{F_0(r_a)}{F_0(\gamma r_b)} \leq \lim_{\eps\to 0} \frac{F_0((\gamma +\eta)r_b)}{F_0(\gamma r_b)} = \left(\frac{\gamma + \eta}{\gamma}\right)^{\alpha}.\]
We used that, since $\psi_r(s)>0$ for $s\in (0,1]$, we have that $r_a\to 0$ and $r_b\to 0$. The inequality holds for all $\eta>0$, and we can show a similar inequality for the $\liminf$, which means that
\[ \lim_{\eps \to 0} \frac{F_0(r_a)}{F_0(\gamma r_b)}  = 1.\]
Since $r_a$ and $r_b$ are decreasing continuous functions of $\eps$, we have shown that in fact
\[ \lim_{t\downarrow 0} \frac{F_0(\gamma t)}{F_0(-t)} = \gamma^{1/2}.\]
This in turn implies that
\[ \lim_{t\downarrow 0} \frac{F_0(-st)}{F_0(-t)} = \lim_{t\downarrow 0} \frac{F_0(+\gamma st)}{F_0(+\gamma t)} = s^\alpha.\]
So the rescaled behavior of $F_0$ to the left of zero is equal to the behavior of $F_0$ to the right of zero.\\

\noindent The rest of the proof is based on Equation \eqref{eq:rewrate}:
\begin{eqnarray*}
\P(\hat{f}(0)\geq a) & = & \P\left( \inf_{-r_a^{-1}\leq s\leq 0}(W_{s} + \eps^{-1}r_a^{-1/2}F_0(r_a s) - Cs) \leq \right.\\
&  & \hspace{2cm} \left. \inf_{0< s\leq r_a^{-1}}(W_{s} + \eps^{-1}r_a^{-1/2}F_0(r_a s) - Cs)\right).
\end{eqnarray*}
Here, $W_s$ is twosided Brownian motion. Note that we can rewrite this equation as
\[ \P(\hat{f}(0)\geq a) = \P \left( \argmin_{s\in [-r_a^{-1}, r_a^{-1}]} (W_{s} + \eps^{-1}r_a^{-1/2}F_0(r_a s) - Cs) \leq 0\right).\]
Using Lemma \ref{lem:limpsi}, we conclude that there exists a family of functions $\eta_t(s)$ on $[0,\infty)$, such that $\eta_t \to 0$ uniformly on compacta as $t\to 0$, with
\begin{equation*}
F_0(st) = s^\alpha F_0(t) + \eta_t(s)F_0(t)\ \ \ \ (t\in \R).
\end{equation*}
This shows that for $s\in [0,\infty)$, we have
\begin{eqnarray*}
W_{s} + \eps^{-1}r_a^{-1/2}F_0(r_a s) - Cs & = & W_s + Cs^\alpha + C\eta_t(s) - Cs\\
& \longrightarrow & W_s + Cs^\alpha -Cs,
\end{eqnarray*}
uniformly on compacta. For $s\in (-\infty, 0]$, we have to be a bit more careful:
\begin{eqnarray*}
\eps^{-1}r_a^{-1/2}F_0(r_a s) & =  & \eps^{-1}r_a^{-1/2}\left(\frac{r_a}{r_b}\right)^\alpha |s|^\alpha F_0(-r_b) + \eps^{-1}r_a^{-1/2}\eta_t(|s|r_a/r_b)F_0(-r_b) \\
& = & \left(\frac{r_a}{r_b}\right)^{\alpha-1/2}|s|^\alpha \eps^{-1}r_b^{1/2}b + \left(\frac{r_a}{r_b}\right)^{-1/2}\eta_t(|s|r_a/r_b)\eps^{-1}r_b^{1/2}b \\
& \to &  C\gamma^{\alpha-1/2}|s|^\alpha,
\end{eqnarray*}
uniformly on compacta. We have shown that uniformly on compacta
\[ W_{s} + \eps^{-1}r_a^{-1/2}F_0(r_a s) - Cs \to \left\{ \begin{array}{ll}
W_s + Cs^\alpha -Cs & \mbox{for } s\geq 0,\\
W_s + C\gamma^{\alpha-1/2}|s|^\alpha - Cs & \mbox{for } s\leq 0, \end{array}\right. \]
Now we wish to use Theorem 2.7 from \mycite{kimpollard:90}, p.198. This Theorem implies that the location of the minimum of the process $W_{s} + \eps^{-1}r_a^{-1/2}F_0(r_a s) - Cs$ converges in distribution to the location of the minimum of its limiting process, provided that this location is $O_p(1)$. To show this last condition, we consider for $M>1$
\[ \P \left(\argmin_{s\in [-r_a^{-1}, r_a^{-1}]}(W_{s} + \eps^{-1}r_a^{-1/2}F_0(r_a s) - Cs) > M\right) \leq \P \left(\inf_{s\geq M} (W_{s} + \eps^{-1}r_a^{-1/2}F_0(r_a s) - Cs) < 0\right).\]
Now we use that for $\eps$ small enough, $F_0(Mr_a)\geq M^\alpha F_0(r_a) - F_0(r_a)$, so for $s\geq M$, using convexity of $F_0$, we get
\[ F_0(sr_a)\geq sF_0(Mr_a)/M \geq M^{\alpha-1}ar_as -ar_as/M.\]
Using this we get
\begin{eqnarray*}
\P \left(\inf_{s\geq M} (W_{s} + \eps^{-1}r_a^{-1/2}F_0(r_a s) - Cs) < 0\right) & \leq & \P \left(\inf_{s\geq M} (W_{s} + CM^{\alpha -1}s - Cs/M - Cs ) < 0\right).
\end{eqnarray*}
Clearly, this last probability goes to zero exponentially fast as $M\to +\infty$, since $\alpha>1$. Now we have to check the lower bound for the location of the minimum:
\begin{eqnarray*}\P \left(\argmin_{s\in [-r_a^{-1}, r_a^{-1}]}(W_{s} + \eps^{-1}r_a^{-1/2}F_0(r_a s) - Cs) <- M\right) & \leq & \P \left(\inf_{s\leq -M} (W_{s} + \eps^{-1}r_a^{-1/2}F_0(r_a s) - Cs) < 0\right)\\
& \leq & \P \left(\inf_{s\leq -M} (W_{s} - Cs) < 0\right).
\end{eqnarray*}
This last probability again goes to zero exponentially fast as $M\to +\infty$. This proves the Theorem for $\gamma > 0$. When $\gamma=0$, so $r_a/r_b\to 0$, the above reasoning goes through, except for the convergence of the process $W_s + \eps^{-1}r_a^{-1/2}F_0(r_a s) - Cs$ for $s\in (-\infty, 0]$. We need to show that
\[ \eps^{-1}r_a^{-1/2}F_0(r_a s) \to 0,\]
uniformly on compact subsets of $(-\infty, 0]$. Fix a compact set $[-M,0]$ and choose $\eps$ so small, that $Mr_a\leq r_b$. Then for all $s\in [-M,0]$,
\begin{eqnarray*}
|\eps^{-1}r_a^{-1/2}F_0(r_a s)| & \leq & \eps^{-1}r_a^{-1/2}F_0(-r_b)|s|\frac{r_a}{r_b}\\
& \leq & C\left(\frac{r_a}{r_b}\right)^{1/2}M\\
& \to & 0.
\end{eqnarray*}

\noindent
Finally we need to prove the last statement. For this, we directly use Equation \eqref{eq:max<0}:
\[ \P(\hat{f}(0)\geq 0) = \P\left( \inf_{-1\leq t\leq 0}(\eps W_t+F_0(t)) \leq \inf_{0< t\leq 1}(\eps W_t+F_0(t))\right).\]
Now we take the usual rescaling, replacing $t$ by $r_as$ and multiplying with $r_a^{-1/2}$:
\begin{eqnarray*}
\P(\hat{f}(0)\geq 0) & = & \P\left( \inf_{-r_a^{-1}\leq s\leq 0}(W_{s} + \eps^{-1}r_a^{-1/2}F_0(r_a s)) \leq \inf_{0< s\leq r_a^{-1}}(W_{s} + \eps^{-1}r_a^{-1/2}F_0(r_a s))\right).
\end{eqnarray*}
Choose $\eps$ small such that $r_a\geq r_b$. Then if $s\leq -r_b/r_a$, we have $F_0(r_as)\geq |s|F_0(-r_b)r_a/r_b$, whereas if $-r_b/r_a\leq s\leq 0$, we still have that $F_0(r_as)\geq 0$, so
\begin{eqnarray*}
\P(\hat{f}(0)\geq 0) & \leq & \P\left( \inf_{-r_a^{-1}\leq s\leq -r_b/r_a}(W_{s} + \left(\frac{r_a}{r_b}\right)^{1/2}C|s|) \leq \inf_{0< s\leq 1}W_{s} + Cs \right)\\
& & +\ \P\left( \inf_{-r_b/r_a\leq s\leq 0}W_{s}  \leq \inf_{0< s\leq 1}W_{s} + Cs\right).
\end{eqnarray*}
Since $r_a/r_b\to +\infty$, these two probabilities clearly go to zero, since $\inf_{0\leq s\leq 1} W_s+Cs < 0$ with probability 1. Note that for this last result, we do not need any other assumptions on $F_0$.\hfill $\Box$\\

\noindent
As before, we introduce the auxiliary function $G_0$ and $H_0$, but now on a full neighborhood of $0$: fix $\delta>0$ and for $t\in (-\delta, \delta)$
\[ G_0(t) = F_0(t)/t\ \ \ \ \mbox{and}\ \ \ \ H_0(t) = t\sqrt{|G_0^{-1}(t)|}.\]
As before, we have that both $G_0$ and $H_0$ are strictly increasing functions on $(-\delta,\delta)$. We also know that the rate equations \eqref{eq:rate} imply that
\[ H_0(a) = C\eps\ \ \mbox{and}\ \ H_0(-b) = -C\eps.\]
\begin{coro}\label{cor:limdist}
Suppose
\[ \lim_{\eps\to 0} \frac{r_a}{r_b} = \gamma,\]
with $\gamma \in [0,\infty)$. Furthermore, suppose that for $s\geq 0$,
\[ \lim_{t\downarrow 0} \frac{F_0(st)}{F_0(t)} = s^{\alpha}\]
for $\alpha> 1$ (see also Lemma \ref{lem:limpsi}). If $W_s$ $(s\in \R)$ denotes twosided standard Brownian motion, define the process
\[ X(s) = \left\{ \begin{array}{ll}
W_s + s^\alpha & \mbox{for } s\geq 0,\\
W_s + \gamma^{\alpha-1/2}|s|^\alpha & \mbox{for } s\leq 0, \end{array}\right.\]
and the process $\hat{X}(s)$ as the greatest convex minorant of $X$.
Then
\[ \eps^{-1}H_0(\hat{f}(0)) \stackrel{d}{\longrightarrow} {\rm sgn}\left(\frac{d\hat{X}}{ds}(0)\right)\left|\frac{d\hat{X}}{ds}(0)\right|^{\frac{2\alpha-1}{2\alpha-2}}.\]
Here, ${\rm sgn}(x)$ denotes the sign of $x\in \R$.
\end{coro}
{\bf Proof:} We start by considering $\P(\eps^{-1}H_0(\hat{f}(0))\geq C)$, for $C>0$. We get
\begin{eqnarray*}
\P(\eps^{-1}H_0(\hat{f}(0))\geq C) & = & \P(\hat{f}(0)\geq a)\\
& \longrightarrow & \P \left( \inf_{s\leq 0}( W_s + C\gamma^{\alpha-1/2} |s|^{\alpha} - Cs ) \leq \inf_{s\geq 0} ( W_s + C |s|^{\alpha} - Cs ) \right),
\end{eqnarray*}
according to Theorem \ref{thm:limdist}. Now replace $s$ by $C^{2/(1-2\alpha)}s$, multiply left and right by $C^{-1/(1-2\alpha)}$ and use Brownian scaling to get
\[ \P(\eps^{-1}H_0(\hat{f}(0))\geq C) \longrightarrow \P \left( \inf_{s\leq 0}( W_s + \gamma^{\alpha-1/2} |s|^{\alpha} - C^{\frac{2\alpha-2}{2\alpha-1}}s ) \leq \inf_{s\geq 0} ( W_s +  |s|^{\alpha} - C^{\frac{2\alpha-2}{2\alpha-1}}s ) \right).\]
Using the switch relation for the greatest convex minorant, we see that
\begin{eqnarray*}
\P(\eps^{-1}H_0(\hat{f}(0))\geq C) & \longrightarrow & \P\left( \frac{d\hat{X}}{ds}(0) \geq C^{\frac{2\alpha-2}{2\alpha-1}}\right)\\
& = & \P\left( {\rm sgn}\left(\frac{d\hat{X}}{ds}(0)\right)\left|\frac{d\hat{X}}{ds}(0)\right|^{\frac{2\alpha-1}{2\alpha-2}} \geq C\right).
\end{eqnarray*}
When $\gamma = 0$, the proof is finished, since in that case
\[ \P\left( \frac{d\hat{X}}{ds}(0) \geq 0\right) = 1.\]
Now suppose $\gamma > 0$. We have seen in the proof of Theorem \ref{thm:limdist} that the scaling of $F_0$ to the left of $0$ is the same as the scaling to the right, so for all $s\geq 0$
\[ \lim_{t\to 0} \frac{F_0(st)}{F_0(t)} = s^\alpha.\]
Consider for $C>0$
\begin{eqnarray*}
\P(\eps^{-1}H_0(\hat{f}(0))\leq -C) & = & \P(\hat{f}(0)\leq -b)\\
& \longrightarrow & \P \left( \inf_{s\leq 0}( W_s + C\gamma^{-\alpha+1/2} |s|^{\alpha} - Cs ) \leq \inf_{s\geq 0} ( W_s + C |s|^{\alpha} - Cs ) \right),
\end{eqnarray*}
using Theorem \ref{thm:limdist} for the left hand side of the origin (that is, interchange $a$ and $b$ and replace $\gamma$ by $1/\gamma$). Now replace $s$ by $-\gamma C^{2/(1-2\alpha)}s$, multiply left and right by $\gamma^{-1/2}C^{-1/(1-2\alpha)}$ and use Brownian scaling to get
\[  \P(\eps^{-1}H_0(\hat{f}(0))\leq -C) \longrightarrow \P \left( \inf_{s\geq 0}( W_s + |s|^{\alpha} + C^{\frac{2\alpha-2}{2\alpha-1}}s ) \leq \inf_{s\leq 0} ( W_s +  \gamma^{\alpha-1/2}|s|^{\alpha} + C^{\frac{2\alpha-2}{2\alpha-1}}s ) \right).\]
Note that the two infima have switched sides because of the scaling with a negative constant. Again using the switch relation we get
\begin{eqnarray*}
\P(\eps^{-1}H_0(\hat{f}(0))\leq -C) & \longrightarrow & \P\left( \frac{d\hat{X}}{ds}(0) \leq -C^{\frac{2\alpha-2}{2\alpha-1}}\right)\\
& = & \P\left( {\rm sgn}\left(\frac{d\hat{X}}{ds}(0)\right)\left|\frac{d\hat{X}}{ds}(0)\right|^{\frac{2\alpha-1}{2\alpha-2}} \leq -C\right).
\end{eqnarray*}
This proves the corollary.\hfill $\Box$\\

\noindent
The condition that $F_0$ is regularly varying around $0$ with parameter $\alpha>1$, implies that the function $H_0$ is regularly varying around $0$ with parameter $\beta=(2\alpha-1)/(2\alpha-2)$, so for all $s\geq 0$
\[ \lim_{t\to 0} \frac{H_0(st)}{H_0(t)} = s^\beta.\]
It is well known from the theory of regularly varying functions that this limit is uniform for $s\in [1/M,M]$, for any $M>1$. This will help us prove the next corollary:
\begin{coro}\label{cor:limdistfin}
With the conditions and notations from Corollary \ref{cor:limdist}, we can show that
\[ \frac{\hat{f}(0)_+}{H_0^{-1}(\eps)} \stackrel{d}{\longrightarrow} \frac{d\hat{X}}{ds}(0)_+ \ \ \mbox{and}\ \  \frac{\hat{f}(0)_-}{-H_0^{-1}(-\eps)} \stackrel{d}{\longrightarrow} \frac{d\hat{X}}{ds}(0)_-.\]
\end{coro}
{\bf Proof:} We wish to show that
\begin{equation}\label{eq:fracprob1}
\frac{H_0(\hat{f}(0))}{\eps}\cdot \left|\frac{H_0^{-1}\left({\rm sgn}(\hat{f}(0))\eps\right)}{\hat{f}(0)}\right|^\beta {\rm sgn}(\hat{f}(0))\to 1\ \ \ \mbox{in probability}.
\end{equation}
Suppose $\eta > 0$. Using Corollary \ref{cor:limdist}, there exists $M>1$ such that for all $\eps$ small enough
\[ \P\left(H_0(\hat{f}(0))\in [-\eps M,-\eps/M]\cup [\eps/M,\eps M]\right) \geq 1-\eta.\]
If $H_0(\hat{f}(0))\in [\eps/M,\eps M]$, we know that
\[ \frac{H_0^{-1}(\eps /M)}{H_0^{-1}(\eps)} \leq \frac{\hat{f}(0)}{H_0^{-1}(\eps)} \leq \frac{H_0^{-1}(\eps M)}{H_0^{-1}(\eps)}.\]
Since $H_0^{-1}$ is regularly varying around $0$ with parameter $1/\beta$, we then know that for $\eps$ small enough,
\[ \frac12 M^{-1/\beta} \leq \frac{\hat{f}(0)}{H_0^{-1}(\eps)} \leq 2M^{1/\beta}.\]
A similar reasoning shows that if $H_0(\hat{f}(0))\in [-\eps/M,-\eps M]$, then for $\eps$ small enough,
\[ \frac12 M^{-1/\beta} \leq \frac{\hat{f}(0)}{H_0^{-1}(-\eps)} \leq 2M^{1/\beta}.\]
Now consider
\[ \eps^{-1}H_0(\hat{f}(0)) = {\rm sgn}(\hat{f}(0)){H_0\left(H_0^{-1}({\rm sgn}(\hat{f}(0))\eps)\,\frac{\hat{f}(0)}{H_0^{-1}({\rm sgn}(\hat{f}(0))\eps)}\right)}/{H_0(H_0^{-1}({\rm sgn}(\hat{f}(0))\eps))}.\]
Since $H_0(st)/H_0(t)\to s^\beta$ uniform for $s$ in compact subsets of $(0,\infty)$, we can conclude with probability higher than $1-\eta$, that for $\eps$ small enough,
\[ \left|\eps^{-1}H_0(\hat{f}(0)) - {\rm sgn}(\hat{f}(0))\left(\frac{\hat{f}(0)}{H_0^{-1}({\rm sgn}(\hat{f}(0))\eps)}\right)^\beta\right|<\eta/M\ \ \mbox{and}\ \ \left|\eps^{-1}H_0(\hat{f}(0))\right|\geq 1/M.\]
This proves \eqref{eq:fracprob1}. Corollary \ref{cor:limdist} then immediately shows that
\[ \frac{{\rm sgn}(\hat{f}(0))\hat{f}(0)}{H_0^{-1}({\rm sgn}(\hat{f}(0))\eps)} \stackrel{d}{\longrightarrow} \frac{d\hat{X}}{ds}(0).\]
This can be written in a nicer way when we look at $\hat{f}(0)_+$ and $\hat{f}(0)_-$:
\[ \frac{\hat{f}(0)_+}{H_0^{-1}(\eps)} \stackrel{d}{\longrightarrow} \frac{d\hat{X}}{ds}(0)_+\]
and
\[ \frac{\hat{f}(0)_-}{-H_0^{-1}(-\eps)} \stackrel{d}{\longrightarrow} \frac{d\hat{X}}{ds}(0)_-.\]
\hfill $\Box$\\

\noindent
Suppose $f_0$ is differentiable in $0$ with $f_0'(0)>0$. Then
\[ \frac{F_0(st)}{F_0(t)} = \frac{\frac12 s^2t^2f'_0(0) + o(t^2)}{\frac12 t^2f'_0(0) + o(t^2)} \to s^2\ \ \ \ (t\to 0).\]
Furthermore, $G_0(t)=F_0(t)/t=\frac12 f'_0(0) t + o(t)$, which implies that $G_0^{-1}(t)=2f'_0(0)^{-1}t + o(t)$, so
\[ H_0(t) = \sqrt{2}f'_0(0)^{-1/2}t^{3/2} + o(t^{3/2}).\]
This means that
\[ H_0^{-1}(\eps) = \left(\frac12 f'_0(0)\right)^{1/3}\eps^{2/3} + o(\eps^{2/3}).\]
Define $X(s)=W_s+s^2$, with $W_s$ twosided Brownian motion, and define $\hat{X}$ as the greatest convex minorant of $X$. Then Corollary \ref{cor:limdistfin} tells us that
\[ \left(\frac12 f'_0(0)\right)^{-1/3}\eps^{-2/3}\hat{f}(0) \stackrel{d}{\longrightarrow} \frac{d\hat{X}}{ds}(0),\]
in accordance with the classical result by Brunk in \mycite{brunk:70}, when translated to the white noise model, except that we do not need a continuous derivative of $f_0$ in a neighbourhood of $0$, we just need the existence of the derivative in $0$.

\subsection{Optimality of the rate}

We wish to show that the rate for the LS-estimator is ``locally optimal'' in the following (non-precise) sense: for each monotone $L^2$-function $f_0$, there exists a sequence of alternative monotone $L^2$-functions $f_1$, such that the rate of the LS-estimator for $f_0$ and $f_1$ cannot both be significantly improved by any other estimator. To be more precise, we will prove the following theorem:
\begin{thm}\label{thm:locrate}
Choose two significance levels $\alpha \in (0,1)$ and $\beta \in (0,1/2)$. There exist $\eta>0$, such that for all $\eps>0$ small enough, we can find a monotone $L^2$-function $f_1$ (close to $f_0$), and we can find a rate $\gamma(\eps)$ with
\[ \limsup_{\eps \to 0}\ \max_{i=0,1}\ \P_{f_i}\left(|\hat{f}(0)-f_i(0)|\geq \gamma(\eps)\right) \leq \alpha\]
and
\[  \liminf_{\eps \to 0}\ \inf_{\hat{\theta}}\ \max_{i=0,1}\ \P_{f_i}\left(|\hat{\theta}(Y)-f_i(0)|\geq \eta\cdot \gamma(\eps)\right) > \beta,\]
where $\hat{\theta}(Y)$ is any estimator of $f(0)$ based on the data $Y$.
\end{thm}
{\bf Remark 1:} One may want to choose different rate-functions $\gamma_0$ and $\gamma_1$ for the two different functions $f_0$ and $f_1$, but it seemed natural to take them equal. In any case, this statement is stronger.\\
{\bf Remark 2:} Choose an event $A\subset C([-1,1])$ such that $\P_{f_0}(Y\in A)\geq 1/2$ and $\P_{f_1}(Y\notin A)\geq 1/2$. Define the estimator
\[ \hat{\theta}(Y) = f_0(0)1_A(Y) + f_1(0)1_{A^c}(Y).\]
Then for any choice of $\eta$ and $\gamma$, we would have
\[ \max_{i=0,1}\ \P_{f_i}\left(|\hat{\theta}(Y)-f_i(0)|\geq \eta\cdot \gamma(\eps)\right) \leq \frac12,\]
which is why in Theorem \ref{thm:locrate}, we choose $\beta \in (0,1/2)$.\\

\noindent {\bf Proof of Theorem \ref{thm:locrate}:} Choose $\eps >0$ small enough such that the equations
\[ F_0(r_a)= ar_a,\ \ F_0(-r_b)= br_b\ \ \mbox{and}\ \ r_a^{1/2}a=r_b^{1/2}b=C\eps\]
have solutions for some fixed $C>0$ with
\[ \frac{2}{\sqrt{2\pi}\,C} \leq \alpha.\]
Suppose that for this $\eps$, $a \geq b$. The case $a < b$ can be handled analogously.
Define for some fixed $0<\delta \leq 1$
\[ f_1(t) = \left\{ \begin{array}{rl}
 \delta a & \mbox{if } t\geq 0\ \mbox{and } f_0(t)\leq \delta a,\\
f_0(t) & \mbox{otherwise}.
\end{array} \right. \]
Then $f_1$ is a monotone $L^2$-function. Note that $f_1$ will be discontinuous in $0$. Define
\[ \gamma(\eps) = 2a.\]
Then Theorem \ref{thm:wnuppbnd} together with Lemma \ref{lem:WpsiC} shows that (remember that $a\geq b$)
\[ \P_{f_0}\left(|\hat{f}(0)|\geq \gamma(\eps)\right) \leq \alpha.\]
Since $f_1\geq f_0$, it easily follows that
\[ \P_{f_1}\left(\hat{f}(0)-\delta a\leq -2a\right) \leq \P_{f_0}\left(\hat{f}(0)\leq -a\right) \leq \frac{1}{\sqrt{2\pi}\,C}.\]
Now we focus on  $\P_{f_1}\left(\hat{f}(0)\geq (2+\delta)a\right)$. Define
\[ F_1(t) = \int_0^t f_1(s)\,ds.\]
We can use Equation \eqref{eq:rewrate} for the situation where the underlying function is $f_1$:
\begin{eqnarray*}
\{\hat{f}(0)\geq 2a\} & = & \{ \inf_{-r_a^{-1}\leq s\leq 0}(r_a^{-1/2}W_{r_a s} + \eps^{-1}r_a^{-1/2}F_1(r_a s) - 2\eps^{-1}r_a^{1/2}as) \leq \nonumber \\
&  & \hspace{2cm} \inf_{0< s\leq r_a^{-1}}(r_a^{-1/2}W_{r_a s} + \eps^{-1}r_a^{-1/2}F_1(r_a s) - 2\eps^{-1}r_a^{1/2}as)\}.
\end{eqnarray*}
Again we have that $F_1(s)\geq 0$ for $s\leq 0$. Define $s_a=\inf\{t>0\ :\ f_0(t)\geq a\}$. Clearly, $s_{\delta a}\leq r_a$. We easily check that for $s\geq s_{\delta a}$, $F_1(s) = \delta as_{\delta a} + F_0(s) -F_0(s_{\delta a})$. This implies that $F_1(r_a)\leq (1+\delta)ar_a\leq 2ar_a$. Since $F_1$ is convex, we conclude that for $0\leq s\leq 1$,
\[ F_1(r_as) \leq 2ar_as.\]
Now we can follow the exact same steps as in the proof of Theorem \ref{thm:wnuppbnd}, starting at Equation \eqref{eq:rewrate}, to conclude that
\[ \P_{f_1}\left(\hat{f}(0)\geq 2a\right)\leq \frac{1}{2\sqrt{2\pi}\,C}.\]
This clearly shows that
\[ \P_{f_1}\left(|\hat{f}(0)-f_1(0)|>\gamma(\eps)\right) = \P_{f_1}\left(\hat{f}(0)>(2+\delta)a\right)+\P_{f_1}\left(\hat{f}(0)<-(2-\delta)a\right) \leq \alpha.\]
So we have shown that our rate $\gamma$ satisfies the first requirement of the theorem.\\

\noindent
Now define $\mu$ as the probability measure on $C([-1,1])$ that corresponds to standard two-sided Brownian motion, and denote with $P_0$ and $P_1$ the measures corresponding to the model with $f_0$ and $f_1$ respectively. It is well known that
\[ \frac{dP_i}{d\mu}(W) = \exp\left(\eps^{-1} \int f_i(t)dW(t) - \frac12 \eps^{-2}\int f_i(t)^2dt\right).\]
Therefore
\[ \frac{dP_1}{dP_0}(W) = \exp\left(\eps^{-1} \int (f_1(t)-f_0(t))dW(t) - \frac12 \eps^{-2}\int f_1(t)^2dt + \frac12 \eps^{-2}\int f_0(t)^2dt\right).\]
This means that
\begin{eqnarray*}
\|P_1-P_0\|_1^2 & \leq & \E_{P_0}\left(\frac{dP_1}{dP_0}(W)-1\right)^2\\
&  = & \E_{P_0}\left(\frac{dP_1}{dP_0}(W)\right)^2 - 1\\
& = & \E_\mu\left(\exp\left(\eps^{-1} \int (2f_1(t)-f_0(t))dW(t) - \eps^{-2}\int f_1(t)^2dt + \frac12\eps^{-2}\int f_0(t)^2dt\right)\right)-1\\
& = & \exp\left(\frac12 \eps^{-2} \int (2f_1(t)-f_0(t))^2dt - \eps^{-2}\int f_1(t)^2dt + \frac12\eps^{-2}\int f_0(t)^2dt\right)-1\\
& = & \exp\left(\eps^{-2} \int (f_1(t)-f_0(t))^2dt\right)-1.
\end{eqnarray*}
We immediately see that
\[ \int (f_1(t)-f_0(t))^2dt \leq \delta^2a^2s_{\delta a} \leq \delta^2a^2r_a\]
so we conclude that
\[ \|P_1-P_0\| \leq \sqrt{\exp(C^{2}\delta^2)-1}.\]
Choose $\delta\in (0,1]$ small enough, such that $\|P_1-P_0\| < 2-4\beta$. Choose $\eta =\delta /4$. Denote with $p_i$ the density of $P_i$ with respect to $\mu$ $(i=0,1)$. We have that for any estimator $\hat{\theta}$
\begin{eqnarray*}
\max_{i=0,1} \P_{f_i}\left(|\hat{\theta}(Y)-f_i(0)|\geq 2\eta a\right) & \geq & \frac12 \sum_{i=0}^1 \P_{f_i}\left(|\hat{\theta}(Y)-f_i(0)|\geq 2\eta a\right)\\
& = & \frac12\ \E_\mu \left( 1_{\{|\hat{\theta}(Y)|\geq 2\eta a\}}p_0(W) + 1_{\{|\hat{\theta}(Y)-\delta a|\geq 2\eta a\}}p_1(W)\right)\\
& \geq & \frac12\ \E_\mu\left(\min(p_0(W),p_1(W))\right)\\
& = & \frac12\ (1-\frac12\|P_1-P_0\|)\\
& > & \beta.
\end{eqnarray*}
This proves the theorem.\hfill $\Box$

\section{The LS-estimator with measurements on a grid}\label{sec:grid}

In this section we wish to show that in the model
\[ Y_i = f_0(x_i) + \eps_i,\]
where $x_i=i/n\ (i=-n,\ldots ,n)$ (so our measurements are taken on a grid) and $\eps_i$ iid, we will get results analogous to the white noise model. The key observation is that when we take measurements on a grid, we can represent the least squares estimator $\hat{f}$ as the derivative of a greatest convex minorant, just as we did with the white noise model. When we define
\[ \hat{g}(t) = \sum_{i=-n}^n Y_i1_{\{x_{i-1}<t\leq x_i\}},\]
with $x_{-n-1}=-1-1/n$, and for $s\in [ -1 - 1/n,1]$
\[ \hat{G}(s) = \int_0^s \hat{g}(t)\,dt,\]
we can define
\[ \hat{F}(t) = \sup\{\phi(t)\ |\ \phi \mbox{ affine and } \forall\ -1-\frac1{n}\leq s\leq 1 : \phi(s)\leq \hat{G}(s)\}.\]
Finally, the least squares estimator is defined as
\[ \hat{f}(t) = \lim_{h\downarrow 0} \frac{F(t)-F(t-h)}{h}.\]
Define $a$ and $b$, depending on $n$, as follows:
\begin{equation}\label{eq:defabgrid}
F_0(r_a)= ar_a,\ \ F_0(-r_b)= br_b\ \ \mbox{and}\ \ r_a^{1/2}a=r_b^{1/2}b=Cn^{-1/2}.
\end{equation}
Here, as before, $C>0$ is some fixed constant. We have the following result:
\begin{thm}\label{thm:gridrate}
With the notations as above, suppose ${\rm Var}(\eps_i)=\sigma^2<+\infty$. Then
\[ \limsup_{n\to \infty} \P(\hat{f}(0)\geq a) \leq \P\left(\inf_{s\leq 0} W_{s}- \frac{C}{\sigma}s \leq  \inf_{0< s\leq 1} W_{s} - \frac{C}{\sigma}(s-\psi_r(s))\right)\]
and
\[ \limsup_{n\to \infty} \P(\hat{f}(0)\leq -b) \leq \P\left(\inf_{s\leq 0} W_{s}- \frac{C}{\sigma}s \leq  \inf_{0< s\leq 1} W_{s} - \frac{C}{\sigma}(s-\psi_l(s))\right).\]
\end{thm}
{\bf Proof:} We start by bounding $\P(\hat{f}(0)\geq a)$; the bound for $\P(\hat{f}(0)\leq -b)$ follows completely analogously. As in \eqref{eq:max<0}, we note that
\[ \{\hat{f}(0)\geq a\} = \{ \inf_{-1-1/n\leq t< 0} \left(\hat{G}(t)-at\right) \leq \inf_{0\leq t\leq 1} \left(\hat{G}(t)-at\right) \}.\]
Define
\[ \tilde{f}_0(t) = \sum_{i=-n}^n f_0(x_i)1_{\{x_{i-1}<t\leq x_i\}}.\]
We now use a similar rescaling as with the white noise model, so $t=r_as$ and multiplying left and right with $n^{1/2}r_a^{-1/2}$. Then $\hat{f}(0)\geq a$ precisely when
\begin{eqnarray}\label{eq:scaledgrid}
\inf_{-r_a^{-1}(1+1/n)\leq s <0}\left( r_a^{-1/2}n^{1/2}\int_0^{r_as}\sum_{i=0}^{-n} \eps_i1_{\{x_{i-1}<t\leq x_i\}}\,dt + r_a^{-1/2}n^{1/2}\int_0^{r_as} \tilde{f}_0(t)\,dt - (r_an)^{1/2}as \right) & \leq \nonumber\\
\inf_{0 \leq s \leq r_a^{-1}}\left( r_a^{-1/2}n^{1/2}\int_0^{r_as}\sum_{i=1}^{n} \eps_i1_{\{x_{i-1}<t\leq x_i\}}\,dt + r_a^{-1/2}n^{1/2}\int_0^{r_as} \tilde{f}_0(t)\,dt - (r_an)^{1/2}as \right).
\end{eqnarray}
Since $f_0$ is increasing, it is not hard to see that for $s\geq 0$
\[ \int_0^{r_as} \tilde{f}_0(t)\,dt \leq F_0(r_as)+n^{-1}f_0(\lceil r_as\rceil).\]
Here, $\lceil r_as\rceil$ signifies the first grid-point bigger than $r_as$. As before, we can use Lemma \ref{lem:unifconvpsi}: for any $\tau\in(0,1)$ there exists a positive continuous function $\eta$ with $F_0(t)=0 \Rightarrow \eta(t)=0$, such that for $s\in [0,\tau]$
\[F_0(r_a s) \leq  F_0(r_a)\psi_r(s) + F_0(r_a)\eta(r_a).\]
Here, $\eta(r_a)\to 0$. This means that for $0\leq s\leq \tau$
\[r_a^{-1/2}n^{1/2}\int_0^{r_as} \tilde{f}_0(t)\,dt \leq  C\psi_r(s) + C\eta(r_a) + r_a^{-1/2}n^{-1/2}f_0(\lceil r_as\rceil).\]
Furthermore, when $t\leq 0$, $\tilde{f}_0(t)\leq 0$. This means that \eqref{eq:scaledgrid} implies
\begin{eqnarray*}
 \inf_{s <0}\left( r_a^{-1/2}n^{1/2}\int_0^{r_as}\sum_{i=0}^{-n} \eps_i1_{\{x_{i-1}<t\leq x_i\}}\,dt - Cs \right) \leq  \hspace{7cm}\\
\inf_{0 \leq s \leq \tau}\left( r_a^{-1/2}n^{1/2}\int_0^{r_as}\sum_{i=1}^{n} \eps_i1_{\{x_{i-1}<t\leq x_i\}}\,dt + C(s-\psi_r(s)) + C\eta(r_a) + \frac{a}{C}\,f_0(\lceil r_as\rceil)\right).
\end{eqnarray*}
Note that the process
\[ s\mapsto r_a^{-1/2}n^{1/2}\int_0^{r_as}\sum_{i=-n}^{n} \eps_i1_{\{x_{i-1}<t\leq x_i\}}\,dt\]
converges to $\sigma W_s$, in the topology of uniform convergence on compacta, where $W_s$ is twosided standard Brownian motion; this is because $nr_a = C^2a^{-2} \to +\infty$. Also, $a\to 0$ and $\eta(r_a)\to 0$ as $n\to \infty$, so we conclude that
\[ \limsup_{n\to \infty} \P\left(\hat{f}(0)\geq a\right) \leq \P\left(\inf_{s<0}\left( W_s - \frac{Cs}{\sigma}\right) \leq \inf_{0\leq s\leq \tau} W_s + C(s-\psi_r(s))\right).\]
Since this holds for any $\tau\in (0,1)$, we have proved the theorem. \hfill $\Box$\\

\subsection{Optimality of the rate}

We wish to prove an analog of Theorem \ref{thm:locrate} for the model with observations on a grid. We need an extra condition on the distribution of $\eps_i$. This makes sense, because suppose that $\eps_i\in {\mathbb Z}$ with probability 1, then it would be very easy to distinguish $f_0$ and $f_1$ if $f_1(0)-f_0(0)\notin {\mathbb Z}$. The condition we need is the following:
\begin{itemize}
\item[{\bf (C1)}] The distribution of $\eps_i$, with ${\rm Var}(\eps_i)=\sigma^2<+\infty$, has a density $\phi$ with respect to the Lebesgue measure, such that there exists $M>0$ with
\[ \forall\ a\in \R:\ \int_{-\infty}^\infty \left(\phi^{1/2}(y-a) - \phi^{1/2}(y)\right)^2\,dy \leq Ma^2.\]
\end{itemize}
This condition would follow from Hellinger differentiability at $0$ of the model $a\mapsto \phi(\cdot-a)$.
\begin{thm}\label{thm:locrategrid}
Suppose Condition (C1) holds. Choose two significance levels $\alpha \in (0,1)$ and $\beta \in (0,1/2)$. There exist $\eta>0$, such that for all $n$ large enough, we can find a monotone increasing function $f_1$ (close to $f_0$), and we can find a rate $\gamma_n$ with
\[ \limsup_{n \to \infty}\ \max_{i=0,1}\ \P_{f_i}\left(|\hat{f}(0)-f_i(0)|\geq \gamma_n\right) \leq \alpha\]
and
\[  \liminf_{n \to \infty}\ \inf_{\hat{\theta}}\ \max_{i=0,1}\ \P_{f_i}\left(|\hat{\theta}(Y)-f_i(0)|\geq \eta\cdot \gamma_n\right) > \beta,\]
where $\hat{\theta}(Y)$ is any estimator of $f(0)$ based on the data $Y$.
\end{thm}
{\bf Proof:} We just follow the steps as in the proof of Theorem \ref{thm:locrate}, so we define $a, r_a, b$ and $r_b$ as in \eqref{eq:defabgrid} with $C>0$ such that
\[ \frac{2\sigma}{\sqrt{2\pi}\,C} \leq \alpha.\]
For fixed $n$, suppose that $a\geq b$. Then we define
\[ \gamma_n = 2a\]
and for some fixed $0<\delta\leq 1$
\[ f_1(t) = \left\{ \begin{array}{rl}
 \delta a & \mbox{if } t\geq 0\ \mbox{and } f_0(t)\leq \delta a,\\
f_0(t) & \mbox{otherwise}.
\end{array} \right. \]
Theorem \ref{thm:gridrate} shows that (remembering that $a\geq b$)
\[ \P_{f_0}\left(|\hat{f}(0)-f_0(0)|\geq \gamma_n\right) \leq \alpha.\]
Since $f_1\geq f_0$, we again have that
\[ \P_{f_1}\left(\hat{f}(0)-\delta a\leq -2a\right) \leq \P_{f_0}\left(\hat{f}(0)\leq -a\right) \leq \frac{1}{\sqrt{2\pi}\,C}.\]
To bound $\P_{f_1}\left(\hat{f}(0)\geq (2+\delta)a\right)$, we follow the proof of Theorem \ref{thm:gridrate}, but with $f_0$ replaced with $f_1$. Define
\[ \tilde{f}_1(t) = \sum_{i=-n}^n f_1(x_i)1_{\{x_{i-1}<t\leq x_i\}}\]
and
\[ F_1(t) = \int_0^t f_1(s)\,ds.\]
Then, in the model using $f_1$, $\hat{f}(0)> 2a$ precisely when
\begin{eqnarray*}
\inf_{-r_a^{-1}(1+1/n)\leq s <0}\left( r_a^{-1/2}n^{1/2}\int_0^{r_as}\sum_{i=0}^{-n} \eps_i1_{\{x_{i-1}<t\leq x_i\}}\,dt + r_a^{-1/2}n^{1/2}\int_0^{r_as} \tilde{f}_1(t)\,dt - (r_an)^{1/2}2as \right) & \leq \nonumber\\
\inf_{0 \leq s \leq r_a^{-1}}\left( r_a^{-1/2}n^{1/2}\int_0^{r_as}\sum_{i=1}^{n} \eps_i1_{\{x_{i-1}<t\leq x_i\}}\,dt + r_a^{-1/2}n^{1/2}\int_0^{r_as} \tilde{f}_1(t)\,dt - (r_an)^{1/2}2as \right).
\end{eqnarray*}
Now, following the steps after Equation \eqref{eq:scaledgrid}, and using the fact that, as in the proof of Theorem \ref{thm:locrate}, for $0\leq s\leq 1$,
\[F_1(r_as)\leq 2as,\]
we conclude that $\gamma_n$ satisfies the first requirement of the theorem.\\

\noindent Now we have to note that in the model with measurements on a grid, the data $Y\in \R^{2n+1}$. As a dominating measure $\mu$ we just take the Lebesgue measure, and we get that the density $p_i$ of the data $Y$, when we are in the model
\[Y_j=f_i(x_j)+\eps_j,\]
is given by
\[ p_i(y) = \prod_{j=-n}^{n}\phi(y_j-f_i(x_j)).\]
Define
\begin{eqnarray*}
\Delta_j & = & \int_{-\infty}^\infty \left( \phi^{1/2}(y-f_1(x_j)) - \phi^{1/2}(y-f_0(x_j))\right)^2\,dy\\
&  = & \int_{-\infty}^\infty \left( \phi^{1/2}(y-(f_1(x_j)-f_0(x_j))) - \phi^{1/2}(y)\right)^2\,dy.
\end{eqnarray*}
Let $H^2(p_0,p_1)$ denote the squared Hellinger distance between $p_0$ and $p_1$. Then it is a standard property of the Hellinger distance that
\begin{eqnarray}\label{eq:hellinger}
H^2(p_0,p_1) & \leq & \sum_{j=-n}^n \Delta_j
\end{eqnarray}
But note that when we define $s_a=\inf\{t>0\ :\ f_0(t)\geq a\}\leq r_a$, we get that $\Delta_j=0$ whenever $j<0$ or $j>ns_{\delta a}$. Furthermore, using Condition {\bf (C1)}, we have that
\begin{equation}\label{eq:Delta}
\Delta_j\leq M\delta^2a^2.
\end{equation}
This shows that
\begin{eqnarray}\label{eq:P1P0}
\|P_1-P_0\|_1 & \leq & 2\sqrt{H^2(p_0,p_1)} \nonumber\\
& \leq & 2\sqrt{(ns_{\delta a}+1)M\delta^2a^2}\\
& \leq & 2\sqrt{4nr_aa^2\delta^2M}\nonumber\\
& = & 4C\delta\sqrt{M}\nonumber.
\end{eqnarray}
It follows that we can choose $\delta>0$ small enough such that $\|P_1-P_0\|_1 < 2-4\beta$. The rest of the proof now follows the proof of Theorem \ref{thm:locrate}.\hfill $\Box$\\

\section{The LS-estimator with measurements on random points}\label{sec:rp}

In this section we consider the model
\[ Y_i = f_0(X_i) + \eps_i,\]
where $X_1,\ldots ,X_n$ is an iid sample in $[-1,1]$ with distribution function $G$, independent of the $\eps_i$'s. We again wish to estimate $f_0(0)$, but our LS-estimator is slightly more complicated now. The idea is to identify the order statistic $X_{(i)}$ with $i$, and calculate the least squares estimator as if the measurements were done on the grid $1,\ldots ,n$. So we define for $0\leq t\leq n$:
\[ h(t) = \sum_{i=1}^n Y_i1_{\{i-1< t \leq i\}}\]
and
\[ H(t) = \int_0^t h(s)\,ds.\]
Furthermore, we define
\[ \tilde{F}(t) =  \sup\{\phi(t)\ |\ \phi \mbox{ affine and } \forall\ 0\leq s\leq n : \phi(s)\leq H(s)\}\]
and
\[ \tilde{f}(t) = \lim_{h\downarrow 0} \frac{\tilde{F}(t)-\tilde{F}(t-h)}{h}.\]
Finally, our estimator of $f_0$ is defined by
\[ \hat{f}(t) = \tilde{f}(m)\ \ \ \mbox{with } X_{(m-1)}<t\leq X_{(m)}.\]
Here, $X_{(m)}$ is the $m^{\rm th}$ order statistic of $X_1,\ldots , X_n$.
In order to control the rate of this estimator, we need some control on how the measurement points behave around $0$. We assume the following condition:
\begin{itemize}
\item[{\bf (C2)}] The distribution $G$ of $X_i$ has a density $g$ with respect to the Lebesgue measure in a neighborhood of $0$, such that $g$ is continuous in $0$ and $g(0)>0$.
\end{itemize}
As before, our rate is defined by
\begin{equation*}
F_0(r_a)= ar_a,\ \ F_0(-r_b)= br_b\ \ \mbox{and}\ \ r_a^{1/2}a=r_b^{1/2}b=Cn^{-1/2},
\end{equation*}
where $C>0$ is some fixed constant.
\begin{thm}\label{thm:randdesrate}
With the notations as above, suppose ${\rm Var}(\eps_i)=\sigma^2<+\infty$ and suppose that (C2) holds. Then
\[ \limsup_{n\to \infty}\ \P(\hat{f}(0)\geq a) \leq \P\left(\inf_{s<0} W_s - \frac{Cg(0)^{1/2}s}{\sigma} \leq \inf_{0\leq s\leq 1} W_s - \frac{Cg(0)^{1/2}}{\sigma}\,(s-\psi_r(s))\right)\]
and
\[ \limsup_{n\to \infty}\ \P(\hat{f}(0)\leq -b) \leq \P\left(\inf_{s<0} W_s - \frac{Cg(0)^{1/2}s}{\sigma} \leq \inf_{0\leq s\leq 1} W_s - \frac{Cg(0)^{1/2}}{\sigma}\,(s-\psi_l(s))\right).\]

\end{thm}
{\bf Proof:} We start by bounding $\P(\hat{f}(0)\geq a)$; the bound for $\P(\hat{f}(0)\leq -b)$ follows completely analogously. Define $m$ such that $X_{(m-1)}<0\leq X_{(m)}$; with probability tending to $1$ we can assume that $1<m<n$ (this follows from Condition (C2)). Note that
\begin{eqnarray*}
 \{\hat{f}(0)\geq a\} & = & \{ \tilde{f}(m)\geq a\} \\
  & = & \{ \inf_{0\leq t\leq m-1} \left(H(t)-at\right) \leq \inf_{m\leq t\leq n} \left(H(t)-at\right) \}\\
  & = & \{ \inf_{0\leq t\leq m-1} \left(H(t)-H(m-1)-a(t-m)\right) \leq \inf_{m\leq t\leq n} \left(H(t)-H(m-1)-a(t-m)\right) \}.
\end{eqnarray*}
We again use a similar rescaling to the one we used for the grid model, namely $t=m+nr_as$ and multiplying left and right with $n^{-1/2}r_a^{-1/2}$. Then $\hat{f}(0)>a$ precisely when
\begin{eqnarray}\label{eq:scaledrand}
\inf_{-m(r_an)^{-1}\leq s \leq -(r_an)^{-1}}\left( r_a^{-1/2}n^{-1/2}\int_{m-1}^{m+nr_as}\sum_{i=1}^{m-1} \eps_i1_{\{i-1<t\leq i\}}\,dt\ + \right. \hspace{3cm}& \nonumber\\
\left.r_a^{-1/2}n^{-1/2}\int_{m-1}^{m+nr_as} \sum_{i=1}^{m-1} f_0(X_{(i)})1_{\{i-1<t\leq i\}}\,dt - (r_an)^{1/2}as \right) & \leq \nonumber\\
\inf_{0 \leq s \leq (n-m)(r_an)^{-1}}\left(r_a^{-1/2}n^{-1/2}\int_{m-1}^{m+nr_as}\sum_{i=m}^{n} \eps_i1_{\{i-1<t\leq i\}}\,dt\ + \right. \hspace{3cm}& \nonumber\\
\left. r_a^{-1/2}n^{-1/2}\int_{m-1}^{m+nr_as} \sum_{i=m}^{n} f_0(X_{(i)})1_{\{i-1<t\leq i\}}\,dt - (r_an)^{1/2}as \right) .
\end{eqnarray}
As before, we have that
\[ s\mapsto r_a^{-1/2}n^{-1/2}\int_{m-1}^{m+nr_as}\sum_{i=1}^{n} \eps_i1_{\{i-1<t\leq i\}}\,dt\]
converges to $\sigma W_s$, with $W_s$ two-sided standard Brownian motion. Also, for $i\leq m-1$, $f_0(X_{(i)})\leq 0$. Finally, suppose that $r_a\geq \eta>0$ for all $n\geq 1$. Then $f_0=0$ on $[0,\eta]$, and it becomes very easy to bound the right-hand side of \eqref{eq:scaledrand} if we limit $s$ to this interval, which would get us the desired result (in this case we would have a parametric rate). Now assume that $r_a\to 0$; then we need that
\[ (n-m)(r_an)^{-1}\to +\infty.\]
This is true with probability $1$, since with probability $1$
\[ m/n \to \P(X_1\leq 0)<1.\]
Furthermore, and most importantly, we need to bound for $0\leq s\leq g(0)$
\[ \int_{m-1}^{m+nr_as} \sum_{i=m}^{n} f_0(X_{(i)})1_{\{i-1<t\leq i\}}\,dt \leq  \sum_{i=m}^{m+\lceil nr_as\rceil} f_0(X_{(i)}).\]
Define $k=\lceil nr_as\rceil+1$ and $D=X_{(m+k)}$. When we condition on $D$, we know that $X_{(m)},\ldots X_{(m+k-1)}$ is an iid sample from $G$ restricted to $[0,D]$. This implies, using Chebyshev,
\[ \P\left(\left.\left|\frac1{k^{1/2}}\sum_{i=m}^{m+k-1} f_0(X_{(i)}) - \frac{k^{1/2}}{G([0,D])}\int_0^D f_0(t)dG(t)\right|>\lambda\ \right|\ D\right) \leq \frac{\int_0^D f_0(t)^2dG(t)}{G([0,D])\lambda^2}\leq f_0(D)^2\lambda^{-2}.\]
It is not hard to see that $D\to 0$ almost surely when $n\to +\infty$, uniformly for $s\in [0,g(0)]$, which proves that
\begin{equation}\label{eq:sumintop1}
r_a^{-1/2}n^{-1/2}\sum_{i=m}^{m+\lceil nr_as\rceil} f_0(X_{(i)}) = \frac{r_a^{1/2}n^{1/2}s}{G([0,D])}\int_0^D f_0(t)dG(t) + o_p(1).
\end{equation}
Now $D$ is the position of the $k+1$-th sample point after $0$, and since $k\to \infty$ and $k/n\to 0$, it is not hard to see, keeping in mind Condition (C2), that
\[ D = \frac{r_as}{g(0)} + O_p(r_a^{1/2}n^{-1/2}).\]
Therefore,
\[ G([0,D]) = r_as(1+o_p(1))\]
and
\[ \frac{r_a^{1/2}n^{1/2}s}{G([0,D])}\left|\int_\frac{r_as}{g(0)}^{\frac{r_as}{g(0)}+O_p(r_a^{1/2}n^{-1/2})} f_0(t)dG(t)\right|\leq  f_0\left(\frac{2r_as}{g(0)}\right)g(0)(1 + o_p(1)) = o_p(1).\]
So \eqref{eq:sumintop1} becomes
\[ r_a^{-1/2}n^{-1/2}\sum_{i=m}^{m+\lceil nr_as\rceil} f_0(X_{(i)}) = r_a^{-1/2}n^{1/2}g(0)F_0\left(\frac{r_as}{g(0)}\right)(1+o_p(1)) + o_p(1).\]
Now we use Lemma \ref{lem:unifconvpsi}: for any $\tau\in (0,1)$, there exists a continuous increasing function $\eta$ on $[0,1]$ with $F_0(t)=0\implies \eta(t)=0$, such that for $0\leq s\leq \tau g(0)$
\[ F_0\left(\frac{r_as}{g(0)}\right) \leq \psi_r(s/g(0))F_0(r_a) + \eta(r_a)F_0(r_a).\]
So finally we conclude that
\begin{eqnarray*}
\limsup_{n\to \infty} \P\left(\hat{f}(0)\geq a\right) &  \leq & \P\left(\inf_{s<0}\left(\sigma W_s - Cs\right) \leq \inf_{0\leq s\leq \tau g(0)} \sigma W_s - C(s-g(0)\psi_r(s/g(0)))\right)\\
& = & \P\left(\inf_{s<0} W_s - \frac{Cg(0)^{1/2}s}{\sigma} \leq \inf_{0\leq s\leq \tau} W_s - \frac{Cg(0)^{1/2}}{\sigma}\,(s-\psi_r(s))\right).
\end{eqnarray*}
Since this holds for any $\tau\in (0,1)$, the theorem follows. \hfill $\Box$\\

\subsection{Optimality of the rate}

We have an analogue to Theorem \ref{thm:locrategrid} for this setting as well:
\begin{thm}\label{thm:locraterandom}
Suppose Conditions (C1) and (C2) hold. Choose two significance levels $\alpha \in (0,1)$ and $\beta \in (0,1/2)$. There exist $\eta>0$, such that for all $n$ large enough, we can find a monotone increasing function $f_1$ (close to $f_0$), and we can find a rate $\gamma_n$ with
\[ \limsup_{n \to \infty}\ \max_{i=0,1}\ \P_{f_i}\left(|\hat{f}(0)-f_i(0)|\geq \gamma_n\right) \leq \alpha\]
and
\[  \liminf_{n \to \infty}\ \inf_{\hat{\theta}}\ \max_{i=0,1}\ \P_{f_i}\left(|\hat{\theta}(Y,X)-f_i(0)|\geq \eta\cdot \gamma_n\right) > \beta,\]
where $\hat{\theta}(Y,X)$ is any estimator of $f(0)$ based on the data $(Y,X)$.
\end{thm}
{\bf Proof:} We can follow the proof of Theorem \ref{thm:locrategrid} (and of Theorem \ref{thm:locrate}), choosing the same alternative function $f_1$, also using the steps in the proof of Theorem \ref{thm:randdesrate} for the alternative $f_1$, right up to the point where we need to bound $\|P_1-P_0\|$. In the random design case, our data consists of $Y$ and $X$, but when we condition on $X$, we can use the inequalities \eqref{eq:hellinger} and \eqref{eq:Delta}, just by replacing $x_j$ by $X_j$. The only difference is that the number $N$ of $X_j$'s in the interval $[0,s_{\delta a}]$ is random. However, we have excellent control on $N$, and by looking at Equation \eqref{eq:P1P0}, we can see that the relevant bound is given by
\[ \E(\sqrt{N}) \leq \sqrt{4r_{a}ng(0)}\ \ \ \mbox{for all }n\mbox{ big enough.}\]
Our conclusion is again that we can choose $\delta>0$ such that $\|P_1-P_0\|_1<2-4\beta$, after which we can follow the proof of Theorem \ref{thm:locrate}.\hfill $\Box$\\

\section{The Grenander estimator for monotone densities}\label{sec:gren}

In this final section we wish to show that our methods also work for the Grenander estimator of a monotone density. Consider a sample $X_1,\ldots ,X_n$ from a monotone decreasing density $f_0$ on $[-1,\infty)$. Assume that $f_0$ is continuous in $0$; we wish to estimate $f_0(0)$. Let $\F_n$ denote the empirical distribution function of the sample $X_1,\ldots ,X_n$. Define
\[ \hat{F}(t) = \inf \{ \phi(t)\ |\ \phi \mbox{ affine and } \forall\ s\geq -1: \phi(s)\geq \F_n(s)\},\]
so $\hat{F}$ is the smallest concave majorant of $\F_n$. The Grenander estimator is now defined as
\[ \hat{f}(t) = \lim _{h\downarrow 0} \frac{\hat{F}(t+h)-\hat{F}(t)}{h}.\]
To find the rate of the Grenander estimator, we define
\begin{equation}\label{eq:F0gren}
F_0(t) = \int_0^t (f_0(0)-f_0(s))\,ds.
\end{equation}
This is a convex function such that $F_0'(0)=0$. Since $f_0$ is decreasing, instead of increasing, when considering the event $\{\hat{f}(0)\geq f_0(0)+a\}$, we have to look to the {\em left}, instead of the right. This results in reversed rate-equations: define $a,b>0$ such that
\begin{equation*}
F_0(r_b)= br_b,\ \ F_0(-r_a)= ar_a\ \ \mbox{and}\ \ r_a^{1/2}a=r_b^{1/2}b=Cn^{-1/2},
\end{equation*}
for some fixed $C>0$. Again we define
\[ \psi_r(s) = \limsup_{t\downarrow 0}\frac{F_0(st)}{F_0(t)} \ \ \mbox{and}\ \ \psi_l(s) = \limsup_{t\uparrow 0}\frac{F_0(st)}{F_0(t)}.\]
We have the following theorem:
\begin{thm}\label{thm:grenupbnd}
With the notations as above, we have that if $r_a\to 0$ and $r_b\to 0$,
\[ \limsup_{n\to \infty}\ \P(\hat{f}(0)\geq a) \leq \P\left(\inf_{s<0} W_s - \frac{C}{\sqrt{f_0(0)}}\,s \leq \inf_{0\leq s\leq 1} W_s - \frac{C}{\sqrt{f_0(0)}}\,(s-\psi_l(s))\right)\]
and
\[ \limsup_{n\to \infty}\ \P(\hat{f}(0)\leq -b) \leq \P\left(\inf_{s<0} W_s - \frac{C}{\sqrt{f_0(0)}}\,s \leq \inf_{0\leq s\leq 1} W_s - \frac{C}{\sqrt{f_0(0)}}\,(s-\psi_r(s))\right).\]
\end{thm}
{\bf Proof:} As before, we will only show how to bound $\P(\hat{f}(0)-f(0)\geq a)$ (in fact, this corresponds to the inequality for $b$ in the other proofs). Define
\[ F(t) = \int_0^t f_0(s)ds\]
and introduce the notation $\F_n(0,t] = \F_n(t)-\F_n(0)$, and likewise $F(0,t]$ (which is in fact equal to $F(t)$). Note that
\begin{eqnarray*}
\{\hat{f}(0)\geq f_0(0)+a\} & = & \{ \inf_{-1\leq t\leq 0}\left(f_0(0)t + at - \F_n(t)\right)\geq \inf_{t\geq 0}\left(f_0(0)t + at - \F_n(t)\right)\}\\
& = & \{ \inf_{-1\leq t\leq 0}\left(f_0(0)t + at - \F_n(0,t]\right)\geq \inf_{t\geq 0}\left(f_0(0)t + at - \F_n(0,t]\right)\}\\
& = & \{ \inf_{-1\leq t\leq 0}\left(F_0(t) + at + F(0,t] - \F_n(0,t]\right)\geq \inf_{t\geq 0}\left(F_0(t) + at + F(0,t] - \F_n(0,t]\right)\}
\end{eqnarray*}
We choose the scaling $t=r_as$ and multiply left and right with $n^{1/2}r_a^{-1/2}$ to get that $\hat{f}(0)\geq f_0(0)+a$ precisely when
\begin{eqnarray}\label{eq:rescalesample}
\inf_{-r_a^{-1}\leq s\leq 0}\left(n^{1/2}r_a^{-1/2}F_0(r_as) + Cs - n^{1/2}r_a^{-1/2}(\F_n(0,r_as]-F(0,r_as])\right) \geq  \hspace{3cm}\nonumber \\
\inf_{s\geq 0}\left(n^{1/2}r_a^{-1/2}F_0(r_as) + Cs - n^{1/2}r_a^{-1/2}(\F_n(0,r_as]-F(0,r_as])\right).
\end{eqnarray}
Again we use Lemma \ref{lem:unifconvpsi}, but now for the function $\psi_l$: for any $\tau\in (0,1)$, there exists a continuous increasing function $\eta$ on $[0,1]$ with $F_0(t)=0\implies \eta(t)=0$, such that for $-\tau\leq s\leq 0$
\[ F_0\left(r_as\right) \leq \psi_l(-s)F_0(-r_a) + \eta(r_a)F_0(-r_a).\]
We conclude that $\hat{f}(0)\geq f_0(0)+a$ implies
\begin{eqnarray*}
\inf_{-\tau\leq s\leq 0}\left( - n^{1/2}r_a^{-1/2}(\F_n(0,r_as]-F(0,r_as]) + C(s+\psi_l(-s))\right) + C\eta(r_a)\geq  \hspace{3cm}\\
\inf_{s\geq 0}\left(Cs - n^{1/2}r_a^{-1/2}(\F_n(0,r_as]-F(0,r_as])\right).
\end{eqnarray*}
What remains is to show that if $r_a\to 0$, the process
\[ Y_n:s\mapsto n^{1/2}r_a^{-1/2}(\F_n(0,r_as]-F(0,r_as]) = n^{-1/2}r_a^{-1/2}\sum_{i=1}^n \left(1_{\{X_i\in (0,r_as]\}} - F(0,r_as]\right)\]
converges in distribution, in the topology of uniform convergence on compacta, to $f_0(0)^{1/2}W_s$, where $W_s$ is two-sided standard Brownian motion. It seems that the classical approach to this problem is the easiest one: the fact that the finite dimensional marginal distributions converge is a relatively straightforward application of the Central Limit Theorem for triangular arrays, since we have written the process as a rescaled sum of independent zero-mean variables; it uses the fact that $F'(0)=f_0(0)$. For tightness of the sequence $Y_n$ it suffices to show that for all $s_1\leq s\leq s_2$ in a compact set, there exists a constant $M>0$ such that
\begin{equation}\label{eq:tightness} \E\left(\left(Y_n(s)-Y_n(s_1)\right)^2\left(Y_n(s_2)-Y_n(s)\right)^2\right) \leq M(s_2-s_1)^2.
\end{equation}
It is not hard to see that the only relevant terms after taking the expectation are
\[ n^{-2}r_a^{-2}\E\left(\Big(1_{\{X_i\in (r_as_1,r_as]\}} - F(r_as_1,r_as]\Big)^2\Big(1_{\{X_j\in (r_as,r_as_2]\}} - F(r_as,r_as_2]\Big)^2\right)\ \ \ \mbox{with}\ i\neq j,\]
of which there are of the order $n^2$. Since $f_0$ is bounded in a neighborhood of $0$, we can find a constant $\tilde{M}>0$ such that for $n$ big enough,
\[ F(r_as_1,r_as_2] \leq \tilde{M}r_a(s_2-s_1).\]
This leads to \eqref{eq:tightness}. We can finally conclude that
\begin{eqnarray*}
 \limsup_{n\to \infty} \P\left(\hat{f}(0)\geq f_0(0)+a\right) & \leq & \P\left(\inf_{-\tau\leq s\leq 0}\left(f_0(0)^{1/2}W_s + C(s+\psi_l(-s))\right)\geq \inf_{s\geq 0}(Cs + f_0(0)^{1/2}W_s)\right)\\
 & = & \P\left(\inf_{s\leq 0}\left(W_s - \frac{C}{\sqrt{f_0(0)}}\,s\right) \leq \inf_{0\leq s\leq \tau}\left(W_s - \frac{C}{\sqrt{f_0(0)}}(s - \psi_l(s))\right)\right).
\end{eqnarray*}
Since this holds for any $\tau \in (0,1)$, we have proved the theorem.\hfill $\Box$\\

When $r_a\to r_0>0$, the process $Y_n(s)$ does not converge to Brownian motion, but to a rescaled Brownian bridge, depending on $F_0$. However, we would still have that when $C\to \infty$, \[\P(\hat{f}(0)-f_0(0)>a)\to 0,\]
so $a$ is still the correct rate (in this case the parametric rate).

\subsection{Optimality of the rate}

In the monotone decreasing density case we also wish to show that the Grenander estimator has the by now familiar optimality property.
\begin{thm}\label{thm:locrategrenander}
Choose two significance levels $\alpha \in (0,1)$ and $\beta \in (0,1/2)$. There exist $\eta>0$, such that for all $n$ large enough, we can find a monotone decreasing density $f_1$ on $[-1,\infty)$ (close to $f_0$), and we can find a rate $\gamma_n$ with
\[ \limsup_{n \to \infty}\ \max_{i=0,1}\ \P_{f_i}\left(|\hat{f}(0)-f_i(0)|\geq \gamma_n\right) \leq \alpha\]
and
\[  \liminf_{n \to \infty}\ \inf_{\hat{\theta}}\ \max_{i=0,1}\ \P_{f_i}\left(|\hat{\theta}(Y)-f_i(0)|\geq \eta\cdot \gamma_n\right) > \beta,\]
where $\hat{\theta}(Y)$ is any estimator of $f(0)$ based on the data $Y$.
\end{thm}
{\bf Proof:} The proof is very similar to the previous ones, but we need to be more careful when choosing the alternative. Choose $n$ large enough such that the equations
\begin{equation*}
F_0(r_b)= br_b,\ \ F_0(-r_a)= ar_a\ \ \mbox{and}\ \ r_a^{1/2}a=r_b^{1/2}b=Cn^{-1/2},
\end{equation*}
have solutions for some fixed $C>0$ with
\[ \frac{2\sqrt{f_0(0)}}{\sqrt{2\pi}\,C} \leq \alpha.\]
Here, $F_0$ is defined in \eqref{eq:F0gren}. Suppose that for this $n$, $a \geq b$. The case $a < b$ can be handled analogously.
Define for some fixed $0<\delta \leq 1$
\[ f_1(t) = \left\{ \begin{array}{rl}
f_0(0) + \delta a & \mbox{if } t\leq 0\ \mbox{and } f_0(t)\leq f_0(0) + \delta a + \eta_a,\\
f_0(t) - \eta_a & \mbox{if } t\leq 0\ \mbox{and } f_0(t)> f_0(0) + \delta a + \eta_a,\\
f_0(t) & \mbox{if } t>0.
\end{array} \right. \]
Then $f_1$ is a monotone decreasing density, if we choose $\eta_a$ such that $\int_{-1}^\infty f_1(t)\,dt = 1$.  This is always possible for $n$ big enough, unless $r_a\to r_0>0$ (i.e., unless $f_0$ is constant on $[-r_0,0]$). However, in this case $\hat{f}(0)$ estimates $f_0(0)$ with a parametric rate (since we consider $a\geq b$), so the conclusions of the theorem will follow. From now on we will assume that $r_a\to 0$. If $b>a$, we only define $f_1(t)=f_0(t)+\eta_b$ for $t\leq 1$ and $f_0(t)\leq f_0(0) - \delta a - \eta_a$; for $t\geq 1$ we would define $f_1(t)=f_0(t)$. Define
$$s_{\delta a}=\inf \{ t>0: f_0(-t) \geq f_0(0)+\delta a\}.$$
We have seen before that $s_{\delta a}\leq r_{\delta a}\leq r_a$. Also,
\[ \int_{-s_{\delta a}}^0 (f_1(t) - f_0(t))\,dt \leq \delta a s_{\delta a}.\]
This gives us an upper bound for $\eta_a$: if $n$ is big enough, such that $f_0(-1/2)>f_0(0) + \delta a + \eta_a$, then
\[\delta a s_{\delta a} \geq \int_{-1}^{-s_{\delta a}} (f_0(t)-f_1(t))\,dt \geq \int_{-1}^{-1/2} (f_0(t) - f_1(t))\,dt = \frac12 \eta_a,\]
so we conclude that for $n$ big enough
\[ \eta_a \leq 2\delta a r_a.\]
Now define
\[ \gamma_n = 2a.\]
Then Theorem \ref{thm:grenupbnd} shows that (remember that $a\geq b$)
\[ \P_{f_0}\left(|\hat{f}(0)-f_0(0)|\geq \gamma_n\right) \leq \alpha.\]
From the way we defined $f_1$, it is clear that we can define $X^{(1)}\sim f_1$ and couple it to $X\sim f_0$, such that $X^{(1)}=X$ if $X\geq 0$, and $X\leq X^{(1)}\leq 0$ otherwise. So if we consider the empirical distribution functions of two samples of $X$ and $X^{(1)}$, call them $\F_n$ and $\F^{(1)}_n$, then $\F_n^{(1)}(t)=\F_n(t)$ if $t\geq 0$, and $\F_n^{(1)}(t)\leq \F_n(t)$ if $t\in [-1,0]$. Now note that
\begin{eqnarray*}
\{\hat{f}^{(1)}(0)\leq f_1(0)-2a\} & = & \{ \inf_{-1\leq t\leq 0}\left(f_1(0)t - 2at - \F_n^{(1)}(t)\right)\leq \inf_{t\geq 0}\left(f_1(0)t - 2at - \F_n^{(1)}(t)\right)\}\\
 & \subset & \{ \inf_{-1\leq t\leq 0}\left(f_1(0)t - 2at - \F_n(t)\right)\leq \inf_{t\geq 0}\left(f_1(0)t - 2at - \F_n(t)\right)\}\\
 & =& \{\hat{f}(0)-f_0(0)\leq \delta a-2a\}.
\end{eqnarray*}
So we get, using that $\delta \leq 1$,
\[ \P_{f_1}\left(\hat{f}(0)-f_1(0)\leq -2a\right) \leq \P_{f_0}\left(\hat{f}(0)\leq -a\right) \leq \frac{\sqrt{f_0(0)}}{\sqrt{2\pi}\,C}.\]
Now we focus on  $\P_{f_1}\left(\hat{f}(0)\geq (2+\delta)a\right)$. Define
\[ F_1(t) = \int_0^t \left(f_1(0)-f_1(s)\right)\,ds\]
and, with a slight abuse of notation,
\[ F^{(1)}(t) = \int_0^t f_1(s)\,ds.\]
We can use Equation \eqref{eq:rescalesample} for the situation where the underlying function is $f_1$, using the coupled sample $X^{(1)}_1,\ldots ,X^{(1)}_n$: $\hat{f}(0)\geq f_1(0)+2a$ precisely when
\begin{eqnarray*}
\inf_{-r_a^{-1}\leq s\leq 0}\left(n^{1/2}r_a^{-1/2}F_1(r_as) + 2Cs - n^{1/2}r_a^{-1/2}(\F_n^{(1)}(0,r_as]-F^{(1)}(0,r_as])\right) \geq  \hspace{3cm}\\
\inf_{s\geq 0}\left(n^{1/2}r_a^{-1/2}F_1(r_as) + 2Cs - n^{1/2}r_a^{-1/2}(\F_n^{(1)}(0,r_as]-F^{(1)}(0,r_as])\right).
\end{eqnarray*}
Note that since $F_1$ is convex, $F_1(r_as)\leq -F_1(-r_a)s$ for $-1\leq s\leq 0$, and that $F_1(-r_a)\leq F_0(-r_a) + r_a\delta a\leq 2ar_a$. Furthermore, for $s\geq 0$, $F_1(r_as)\geq 0$. This means that $\hat{f}(0)\geq f_1(0)+2a$ implies that
\[ \inf_{-1\leq s\leq 0}\left(- n^{1/2}r_a^{-1/2}(\F_n^{(1)}(0,r_as]-F^{(1)}(0,r_as])\right) \geq  \hspace{3cm}\\
\inf_{s\geq 0}\left(Cs - n^{1/2}r_a^{-1/2}(\F_n^{(1)}(0,r_as]-F^{(1)}(0,r_as])\right).\]
Since the left-derivative of ${F^{(1)}}$ in $0$ equals $f_0(0)+o_p(1)$, we can proceed as in the proof of Theorem \ref{thm:grenupbnd} to conclude that our rate $\gamma_n$ satisfies the first requirement of the theorem.\\

\noindent
Now we need to bound $\|P_1-P_0\|_1$, where $P_0$ is the distribution of $(X_1,\ldots, X_n)$ and $P_1$ the distribution of $(X^{(1)}_1,\ldots, X^{(1)}_n)$. We do this by bounding the Hellinger distance between $f_0$ and $f_1$:
\begin{eqnarray*}
H^2(f_0,f_1) & = & \int_{-1}^0 \left(\sqrt{f_1(s)}-\sqrt{f_0(s)}\right)^2\,ds\\
& = & f_0(0)\,\int_{-1}^0 \left(\sqrt{\frac{f_1(s)}{f_0(0)}}-\sqrt{\frac{f_0(s)}{f_0(0)}}\,\right)^2\,ds\\
& \leq & \frac12 f_0(0) \,\int_{-1}^0 \frac{(f_1(s)-f_0(s))^2}{f_0(0)^2}\,ds\\
& \leq & \frac{a^2\delta^2s_{\delta a} + \eta_a^2}{2f_0(0)}\\
& \leq & \frac{a^2\delta^2r_{a}}{f_0(0)}.
\end{eqnarray*}
For the last inequality we use that for $n$ big enough, $r_a<1/2$. In the case where $b>a$, you could use the fact that for $n$ big enough, $f_0(r_b)\geq f_0(0)/2$, to get the first inequality (with a different constant). It now follows that
\[ \|P_1-P_0\|_1 \leq 2\sqrt{H^2(p_1,p_0)} \leq 2\sqrt{nH^2(f_0,f_1)}\leq C\delta/\sqrt{f_0(0)}.\]
Choose $\delta\in (0,1]$ small enough, such that $\|P_1-P_0\| < 2-4\beta$, and follow the proof of Theorem \ref{thm:locrate}. \mbox{} \hfill $\Box$

\noindent Eric Cator\\
Delft University of Technology\\
Mekelweg 4\\
2628 CD Delft\\
The Netherlands\\
email: e.a.cator@tudelft.nl

\end{document}